\documentclass[english,11pt,a4paper]{article}

\usepackage[english]{babel}
\usepackage{enumitem} 
\usepackage{wrapfig}
\usepackage[final]{graphicx}
\usepackage{epsfig}
\usepackage{textcomp}
\usepackage{amsmath}
\usepackage{amssymb}
\usepackage{amsfonts}
\usepackage{mathbbol} 
\usepackage{psfrag}
\usepackage{latexsym}
\usepackage{eurosym}
\usepackage{color}    
\usepackage{mathtools}
\usepackage{tikz}
\usepackage{stmaryrd} 
\usepackage[normalem]{ulem}
\usepackage[numbers,sort&compress]{natbib}

\newcommand{\epsi}{\varepsilon}
\newcommand{\IC}{\mathbb{C}} 
\newcommand{\IN}{\mathbb{N}}

\newcommand{\IR}{\mathbb{R}}

\newcommand{\tend}{t_{\mbox{\tiny end}}}
\newcommand{\tstar}{t_\star}

\newcommand{\diag}{\mbox{diag}}

\newcommand{\curl}{\operatorname{curl}}
\renewcommand{\div}{\operatorname{div}}

\newcommand{\hs}[1]{\hspace*{#1 mm}}

\newcommand{\Ord}[1]{\mathcal{O} \! \left(#1\right)}

\newcommand{\A}{\mathcal{A}}
\newcommand{\B}{\mathcal{B}}
\newcommand{\F}{\mathcal{F}}
\newcommand{\G}{\mathcal{G}}

\newcommand{\J}{\mathcal{J}}
\renewcommand{\L}{\mathcal{L}}
\renewcommand{\P}{\mathcal{P}}
\renewcommand{\S}{\mathcal{S}}
\newcommand{\T}{\mathcal{T}}
\newcommand{\U}{\mathcal{U}}
\newcommand{\X}{\mathcal{X}}
\newcommand{\TT}{S}
\newcommand{\FF}{\mathbf{F}}

\newcommand{\ii}{\mathrm{i}}

\newcommand{\w}{\omega}
\newcommand{\Chat}{\widehat{C}}
\newcommand{\uhat}{\widehat{u}}
\newcommand{\vhat}{\widehat{v}}
\newcommand{\phat}{\widehat{p}}

\newcommand{\uu}{\mathbf{u}}

\newcommand{\uutilde}{\mathbf{\widetilde{u}}}
\newcommand{\Pperp}{P^\perp}

\newcommand{\pt}{\partial_{t}}

\newcommand{\intl}{\int\limits}

\newcommand{\weg}[1]{}

\newcommand{\mynote}[1]{} 
\newcommand{\myfootnote}[1]{} 
\newcommand{\proof}{\noindent\textbf{Proof. }}
\newcommand{\proofof}[1]{\noindent\textbf{Proof of #1. }}
\def\qed {\hfill \rule{0,2cm}{0,2cm} \bigskip \\}

\newcommand{\il}{{\interleave}}

\newtheorem{Theorem}{Theorem}[section]
\newtheorem{Proposition}[Theorem]{Proposition}
\newtheorem{Lemma}[Theorem]{Lemma}

\newtheorem{Assumption}[Theorem]{Assumption}

\newtheorem{Remark}[Theorem]{Remark}

\newcounter{fig}

\begin{document}

\title{Approximation of high-frequency wave propagation in dispersive media\thanks{Funded by the Deutsche Forschungsgemeinschaft (DFG, German Research Foundation) -- Project--ID 258734477 -- SFB 1173.}}
\author{Julian Baumstark\thanks{Karlsruher Institut f\"ur Technologie, Fakult\"at f\"ur Mathematik, Institut f\"ur Angewandte und Numerische Mathematik, Englerstr.~2, D-76131 Karlsruhe,
\texttt{julian.baumstark@kit.edu}, \texttt{tobias.jahnke@kit.edu}} \and Tobias Jahnke\footnotemark[2]}
\date{\today}
\maketitle

\abstract{
We consider semilinear hyperbolic systems with a trilinear nonlinearity. Both the differential equation and the initial data contain the inverse of a small parameter $\epsi$, and typical solutions oscillate with frequency proportional to $1/\epsi$ in time and space. Moreover, solutions have to be computed on time intervals of length $1/\epsi$ in order to study 
nonlinear and diffractive effects. As a consequence, direct numerical simulations are extremely costly or even impossible.
We propose an analytical approximation and prove that it approximates the exact solution up to an error of $\Ord{\epsi^2}$
on time intervals of length $1/\epsi$. This is a significant improvement over the classical nonlinear Schr\"odinger approximation, which only yields an accuracy of $\Ord{\epsi}$.
}

\paragraph{Keywords:} High-frequency wave propagation, semilinear wave equation, Maxwell--Lorentz system, diffractive geometric optics, slowly varying envelope approximation, error bounds



\section{Introduction}

We consider semilinear hyperbolic systems of the form
\begin{subequations}
\label{PDE.uu}
\begin{align}
\label{PDE.uu.a}
 \pt \uu +A(\partial)\uu +\frac{1}{\epsi} E\uu &= \epsi T(\uu,\uu,\uu),
 & & t\in(0,\tend/\epsi], \; x\in\IR^d,
 \\
\label{PDE.uu.b}
 \uu(0,x) &= p(x)e^{\ii (\kappa \cdot x)/\epsi} + c.c.
\end{align}
\end{subequations}
with a small parameter $0<\epsi \ll 1$ and vector-valued solutions $\uu: [0,\tend/\epsi] \times \IR^d \rightarrow \IR^n$ for some $\tend>0$ and $d,n\in\IN$. The differential operator 
\begin{align}
\label{Def.A.partial}
A(\partial) = \sum_{\ell=1}^d A_\ell \partial_\ell
\end{align}
contains constant symmetric matrices $A_1, \ldots, A_d \in \IR^{n \times n}$. $ E\in\IR^{n\times n} $ is a skew-symmetric matrix and $T: \IR^n \times \IR^n \times \IR^n \rightarrow \IR^n $ is a trilinear nonlinearity. Its trilinear extension to $\IC^n \times \IC^n \times \IC^n$, which will appear later, is also denoted by $T$.
The initial data \eqref{PDE.uu.b} depend on a smooth envelope function $p: \IR^d \rightarrow \IR^n$ and a wave vector $\kappa \in \IR^d\setminus\{0\}$. 
As usual, ``c.c.'' means complex conjugation of the previous term.

A prominent example in this problem class is the Maxwell--Lorentz system
\begin{align*}
 \pt \mathbf{B} &= -\curl \mathbf{E}, \\
 \pt \mathbf{E} &= \curl \mathbf{B} - 
\frac{1}{\epsi}\mathbf{Q}, \\
 \pt \mathbf{Q} &= \frac{1}{\epsi}(\mathbf{E}-\mathbf{P}) + 
\epsi|\mathbf{P}|_2^2 \mathbf{P}, \\
 \pt \mathbf{P} &= \frac{1}{\epsi}\mathbf{Q}, \\
 \div (\mathbf{E}+\mathbf{P}) &= \div \mathbf{B} = 0, 
\end{align*}
\mynote{Die Divergenzbedingungen stehen so im Paper von Joly, Metivier, Rauch: ``Global solvability ...'' und im Schneider-Buch (S. 441). Im abstrakten Problem \eqref{PDE.uu} gibt es keine Divergenzbedingungen. Darauf w\"urde ich hier aber nicht eingehen.}
which models the propagation of a light beam in a Kerr medium; 
cf.~\cite{donnat-rauch:97a,colin-lannes:09,lannes:11,lannes:98,donnat-rauch:97b,joly-metivier-rauch:96,colin-gallica-laurioux:05}. 
As usual, $ \mathbf{E} $ 
and $ \mathbf{B} $ describe the electric and magnetic fields, respectively. 
$\mathbf{P}$ is the polarization and $ \mathbf{Q}/\epsi$ its time derivative.
In this model, the Maxwell equations for $\mathbf{E}$ and $\mathbf{B}$ are coupled to ordinary differential equations for $ \mathbf{P}$ and $\mathbf{Q}$.\myfootnote{Mit $|\mathbf{P}|$ scheint eher die Euklidsche Vektornorm zu sein (vgl. Colin-Lannes S. 713). 
Roland Schnaubelt hat das best\"atigt.}
The equations are normalized in such a way that the speed of light is 1, and
the parameter $\epsi$ corresponds to the ratio between the wavelength of light and the next characteristic length of the problem; cf.~\cite{donnat-rauch:97a}.
The abstract problem setting applies also to the Klein--Gordon system
\begin{align*}
 \pt \uu + \begin{pmatrix} 0 & \nabla \\ \nabla^T & 0 \end{pmatrix}\uu
+ \frac{1}{\epsi} \begin{pmatrix} 0 & -\nu^T \\ \nu & 0 \end{pmatrix}\uu
= \epsi |\uu|_2^2 M\uu
\end{align*}
with $\nu \in \IR^d\setminus \{0\}$ and with a skew-symmetric matrix $M\in\IR^{n \times n}$; cf.~\cite{colin-lannes:09,lannes:11}.
\mynote{(Braucht man Schiefsymmetrie? Nicht ganz klar. 
Zumindest scheint die $L^2$-Norm (nicht $L^1$!) der L\"osung im Fall einer schiefsymmetrischen Matrix erhalten zu bleiben. Bildet man auf beiden Seiten des Klein--Gordon-Systems das Skalarprodukt mit $\uu$, so erh\"alt man
\begin{align*}
 \uu^T \pt \uu + \uu^T \begin{pmatrix} 0 & \nabla \\ \nabla^T & 0 \end{pmatrix}\uu 
 + \frac{1}{\epsi} \underbrace{\uu^T \begin{pmatrix} 0 & -\nu^T \\ \nu & 0 \end{pmatrix}\uu}_{=0}
= \epsi |\uu|_2^2 \underbrace{\uu^T M\uu}_{=0}
\end{align*}
und damit
\begin{align*}
\tfrac{1}{2} \pt \| \uu(t) \|_{L^2}^2
&=
\tfrac{1}{2} \intl_{\IR^d} \pt | \uu(t,x) |_2^2 \; dx
\\
&=
\intl_{\IR^d} \uu^T(t,x) \pt \uu(t,x) \; dx
\\
&=
- \intl_{\IR^d} \uu^T \begin{pmatrix} 0 & \nabla \\ \nabla^T & 0 \end{pmatrix}\uu  \; dx
= 0,
\end{align*}
wobei der letzte Schritt aus partieller Integration folgt.
Allerdings ist das bei Maxwell--Lorentz nicht so, denn dort geht ja nur die Norm $|P|$ in die rechte Seite ein.)
} 

Physically relevant solutions of \eqref{PDE.uu} oscillate rapidly in time and space due to the small parameter $\epsi$ which occurs both in the PDE \eqref{PDE.uu.a} and in the initial data \eqref{PDE.uu.b}.
Moreover, the problem is scaled in such a way that nonlinear and diffractive effects appear only on a long time interval $[0,\tend/\epsi]$.
As a consequence, approximating the solution of \eqref{PDE.uu} numerically with standard methods is prohibitively inefficient or even unfeasible.
This problem has motivated many attempts to devise simpler models which are more suitable for numerical computations and at the same time provide a reasonable approximation to $\uu$. Among these models, the nonlinear Schr\"odinger approximation
is particularly appealing; cf.~\cite{colin:02,colin-lannes:09,lannes:11,donnat-joly-metivier-rauch:96,joly-metivier-rauch:98,kirrmann-schneider-mielke:92,schneider-uecker:17}. 
In this model the exact solution of \eqref{PDE.uu} is approximated by
\begin{align*}
\uu_{\text{\tiny NLS}}(t,x) = 
e^{\ii (\kappa \cdot x - \w t)/\epsi} U_{\text{\tiny NLS}}(t,x) + c.c.,
\end{align*}
where $(\w,\kappa)$ satisfy the dispersion relation and $ U_{\text{\tiny NLS}} $ evolves according to a nonlinear Schr\"o\-dinger equation. 
If a co-moving coordinate system is used, then this PDE does \emph{not} depend on $\epsi$ and has only to be solved on a time interval of length $\tend$ instead of $\tend/\epsi$; cf.~Remark 8.v. in \cite{colin-lannes:09}. In \cite[Corollary 2]{colin-lannes:09} the error bound
\begin{align}
\label{error.bound.nls}
\sup_{t\in[0,\tend/\epsi]} \| \uu(t) - \uu_{\text{\tiny NLS}}(t) \|_{L^\infty(\IR^n)} \leq C\epsi
\end{align}
was shown under a number of assumptions.
Hence, the nonlinear Schr\"odinger approximation offers a possibility to approximate the solution of \eqref{PDE.uu} up to $\Ord{\epsi}$ without difficulties caused by oscillations or a long time interval.
In some situations, however, a more accurate approximation to $\uu$ is required. 
Our goal is to derive a system of PDEs which is numerically more favorable than
\eqref{PDE.uu} but provides an approximation to 
the solution $\uu$ up to an error of  $\Ord{\epsi^2}$. 
In some respects, our ansatz can be considered as a higher-order extension of the classical slowly varying envelope approximation.

Asymptotic expansions of solutions to systems similar to \eqref{PDE.uu} have been derived, e.g., in 
\cite{rauch-book:12,donnat-rauch:97b,joly-metivier-rauch:00,joly-metivier-rauch:93} for geometric optics, i.e. for times of length $\Ord{1}$. 
In contrast to these works, we seek approximations on time intervals of length
$\Ord{1/\epsi}$, which is the regime of diffractive geometric optics.
In the diffractive regime approximations with infinitely small residual have been constructed in \cite{donnat-joly-metivier-rauch:96} for semilinear and quasilinear systems, but with $\epsi E$ instead of $E/\epsi$ in \eqref{PDE.uu.a}.
More general nonlinear hyperbolic systems, but with $E=0$ have been analyzed in \cite{joly-metivier-rauch:98}. 
Approximate solutions for quasilinear systems with dispersion have been analyzed in \cite{lannes:98},
and for dispersive problems with bilinear nonlinearity in \cite{colin:02}, but without an explicit convergence rate.
The approximation of PDEs by nonlinear Schr\"odinger equations and other modulation equations is extensively discussed
in \cite{schneider-uecker:17} and references therein.

It is well-known that the accuracy of the nonlinear Schr\"odinger approximation deteriorates in the case of short or chirped pulses, and for such situations many improved models have been proposed and analyzed, e.g., in \cite{barrailh-lannes:02,colin-lannes:09,lannes:11,colin-gallica-laurioux:05,alterman-rauch:00,alterman-rauch:03,chung-jones-schaefer-wayne:05}. However, this is not the situation we consider here. 
We restrict ourselves to wave trains where the envelope $p$ varies on a scale which is much larger than the wavelength of the oscillations, but we strive for higher accuracy.
\mynote{Das Paper \cite{joly-metivier-rauch:93} gibt es nur als Hardcopy in der Mathebib. Ich habe es mir angesehen und die ersten paar Seiten fotografiert.}

In the next section an approximation to the solution of \eqref{PDE.uu} is constructed; cf.~\eqref{Ansatz}-\eqref{PDE.mfe.initial.data}. Moreover, we formulate a number of assumptions, prove well-posedness of the approximation, and we compile the toolbox used in the analysis in later sections.
Our main result is Theorem~\ref{Theorem.error.bound}, which provides an error bound for the approximation \eqref{Ansatz}. The proof relies on another important result, namely the fact that certain ``parts'' of the constructed approximation are of $\Ord{\epsi}$, roughly speaking. 
This statement is made precise in Propositions~\ref{Proposition01} and \ref{Proposition02}, and Section~\ref{Sec:Refined.bounds} 
is mostly devoted to the proof of these propositions.
A possible extension to higher accuracy is briefly discussed in the last section.

\paragraph{Notation.}
The Euclidean scalar product of vectors $v,w \in \IC^n$ is denoted by $v \cdot w = v^* w $,
and $ |v|_q $ is the usual $q-$norm of $v$.
For functions $f=f(t,x)$ depending on $t$ and $x$ we will often omit the spatial variable and write $f(t)$ instead of $f(t,x)$.
In the same spirit, the second argument of the Fourier transform $\widehat{f}(t,k)$ of such a function will most often be omitted.
The imaginary unit is denoted by $\ii$, whereas $i$ is used as an index in a few formulas.

\section{Analytical setting}\label{Sec.analytical.setting}
\subsection{Assumptions, ansatz and main goal}\label{Subsec.Ansatz}

For $\alpha \in \IR$ and $\beta \in \IR^d$ we define the matrices
\begin{align}
\label{Def.Ak}
 A(\beta)&=\sum_{\ell=1}^d \beta_\ell A_\ell,
 \\
\label{Def.L}
\L(\alpha,\beta) &= -\alpha I + A(\beta)-\ii E \in \IC^{n\times n}.
\end{align}
The notation \eqref{Def.Ak} is consistent with the definition of $A(\partial)$ in \eqref{Def.A.partial}.
$A(\beta)$ is symmetric by definition, and hence $\L(\alpha,\beta)$ is Hermitian.

From now on let $\kappa \in \IR^d\setminus\{0\}$ be a fixed wave vector and let $\w=\w(\kappa)$ be an eigenvalue of $A(\kappa)-\ii E$. Hence, the matrix $\L(\w,\kappa)$ has a non-trivial kernel. The following assumptions are made.

\begin{Assumption}\label{Ass:polarization}
\text{ }
\begin{enumerate}[label=(\roman*)]
\item \label{Ass:polarization.item.i} 
The kernel of $\L(\w,\kappa)$ is one-dimensional. 

\item \label{Ass:polarization.item.ii} 
The initial data \eqref{PDE.uu.b} are polarized, i.e. 
\begin{align*}
 p(x)  \in \text{ker}\big(\L(\w,\kappa)\big) \text{ for all } x \in \IR^d.
\end{align*}

\item \label{Ass:polarization.item.iii} 
For $j\in\{3,5\}$ the matrix $ \L(j\w, j\kappa)$ is invertible. 
\end{enumerate}
\end{Assumption}
Assumption \ref{Ass:polarization.item.i} could be dropped at the cost of a more complicated notation.
Instead of \ref{Ass:polarization.item.ii} it is actually sufficient to assume that $ p(x)=p_0(x)+\epsi p_1(x) $ with 
$p_0(x) \in \text{ker}\big(\L(\w,\kappa)\big) $. 
This assumption has also been made, e.g., in \cite{colin-lannes:09,lannes:11}. 
The assumption that $p_1 = 0$ is made in order to simplify the presentation.
For $j=3$ Assumption \ref{Ass:polarization.item.iii} was also made in \cite[Assumption 3]{colin-lannes:09}.

The matrices \eqref{Def.L} will play an important role throughout. The analysis in later sections requires the following smoothness properties of eigenvalues and eigenvectors.
\begin{Assumption}\label{Ass:L.properties}
\text{ }
\begin{enumerate}[label=(\roman*)]
\item \label{Ass:L.properties.item.i}
The matrix $\L(0,\beta) = A(\beta)-\ii E$ has a smooth eigendecomposition:
if $\lambda(\beta)$ is an eigenvalue of $\L(0,\beta)$, then $ \lambda \in C^\infty(\IR^d\setminus\{0\},\IR) $,
and there is a corresponding eigenvector $ \psi(\beta) $ such that $ |\psi(\beta)|_2=1 $ for all $\beta$ and 
$ \psi \in C^\infty(\IR^d\setminus\{0\},\IC^{n}). $
\mynote{Eigentlich brauchen wir das nur beim CL-Resultat \"uber die Projektion der Anfangswerte.}
\item \label{Ass:L.properties.item.ii}
Every eigenvalue $\lambda(\beta)$ of $\L(0,\beta)$ is globally Lipschitz continuous, i.e. there is a constant $C$ such that
\begin{align*}
| \lambda(\widetilde{\beta}) - \lambda(\beta) | \leq C | \widetilde{\beta} - \beta |_1 \qquad 
\text{for all } \; \widetilde{\beta}, \beta \in\IR^d.
\end{align*}
\end{enumerate}
\end{Assumption}
Assumption \ref{Ass:L.properties.item.i} corresponds to Assumption 2 in \cite{colin-lannes:09}. 
For the Maxwell--Lorentz system and the Klein--Gordon system 
the eigenvalues are stated in \cite[Example 3 and 4]{colin-lannes:09}, and it can be checked that these eigenvalues do indeed have the properties \ref{Ass:L.properties.item.i} and \ref{Ass:L.properties.item.ii}. 
\mynote{Theorem 3.I.1 in \cite{Rau12} liefert Glattheit der Eigenprojektoren, aber nicht der Eigenvektoren. Da Eigenvektoren nur bis auf Vorfaktoren bestimmt sind, k\"onte man die immer unstetig w\"ahlen.}

Note that $ \L(\alpha,\beta) = -\alpha I + \L(0,\beta) $ has the same eigenvectors as $\L(0,\beta)$, and that the eigenvectors are shifted by $-\alpha$. Hence, if Assumption~\ref{Ass:L.properties} is fulfilled, 
then for every $\alpha\in\IR$ the eigenvalues and eigenvectors of $ \L(\alpha,\beta) $ have the smoothness specified in
\ref{Ass:L.properties.item.i} and \ref{Ass:L.properties.item.ii}, too.

\bigskip

In \cite{colin-lannes:09} the classical nonlinear Schr\"odinger approximation is derived in two steps. The first step, known as the \emph{slowly varying envelope approximation}, is to approximate
\begin{align*}
\uu(t,x) \approx 
\uu_{\text{\tiny SVEA}}(t,x) = 
e^{\ii (\kappa \cdot x - \w t)/\epsi} U_{\text{\tiny SVEA}}(t,x) + c.c.,
\end{align*}
where $U_{\text{\tiny SVEA}}: [0,\tend/\epsi] \times \IR^d \rightarrow \IC^n$ is the (complex-valued) solution of a PDE called the envelope equation.
The accuracy of this approximation is $\Ord{\epsi}$ on $[0,\tend/\epsi]$ in suitable norms; cf.~\cite[Theorem 1]{colin-lannes:09}.
In the second step, it is shown that the envelope equation can be replaced by the nonlinear Schr\"odinger equation without spoiling the accuracy.
From this procedure it is clear that in order to achieve a higher accuracy, the slowly varying envelope approximation has to be replaced by a better one.
In fact, substituting $\uu_{\text{\tiny SVEA}}$ into the nonlinearity
$T$ yields with $U=U\text{\tiny SVEA}(t,x)$
\begin{align*}
T(\uu_{\text{\tiny SVEA}},\uu_{\text{\tiny SVEA}},\uu_{\text{\tiny SVEA}})(t,x)
=&
e^{3 \ii (\kappa \cdot x - \w t)/\epsi}
T(U,U,U)
\\
&+
e^{\ii (\kappa \cdot x - \w t)/\epsi}
\Big(T(U,U,\overline{U})+T(U,\overline{U},U)+T(\overline{U},U,U)\Big)
\\
&+
e^{-\ii (\kappa \cdot x - \w t)/\epsi}
\Big(T(U,\overline{U},\overline{U})+T(\overline{U},U,\overline{U},U)+T(\overline{U},\overline{U},U)\Big)
\\
&+ e^{-3 \ii (\kappa \cdot x - \w t)/\epsi}
T(\overline{U},\overline{U},\overline{U}),
\end{align*}
and the terms with prefactor $e^{\pm3 \ii (\kappa \cdot x - \w t)/\epsi}$
(``higher harmonics'') are ignored in the envelope equation.
This is the motivation to make the ansatz
\begin{align}
\label{Ansatz}
 \uu(t,x) &\approx \uutilde(t,x) = 
 \sum_{j\in \J}
e^{\ii j(\kappa \cdot x - \w t)/\epsi} u_j(t,x),
\qquad u_{-j} = \overline{u_j}
\end{align}
for $ \J = \{\pm 1, \pm 3\}$. 
Terms of the form $ e^{2\ii (\kappa \cdot x - \w t)/\epsi} u_2(t,x) $ are not required because $T$ is trilinear.

\mynote{Unterschiede zu dem, was andere machen: 1. Zus\"atzliche Variable $\theta$. 2. Vor jedem Term steht eine $\epsi$-Potenz, und man fordert, dass im Residuum die Terme mit den niedrigsten $\epsi$-Potenzen verschwinden. Wenn wir tats\"achlich erkl\"aren wollen, was andere Leute anders machen, dann m\"ussen wir weiter ausholen. Das sollte man vielleicht eher in der Einleitung machen.}
\mynote{Usually the ansatz is substituted into \eqref{PDE.uu}, which yields a formal expansion of the residual in powers of $\epsi$, and equations for $u_1$ and $u_3$ are then obtained by imposing the condition that the leading order terms of the residual vanish. We proceed in a different way.}
\mynote{Manche Leute (z.B. Colin-Lannes) machen das so und bekommen trotzdem eine Approximation (z.B. NLS) ohne algebraische Nebenbedingungen, d.h. f\"ur beliebige Anfangsdaten. Deshalb w\"urde ich das hier nicht thematisieren.
}

Substituting the ansatz \eqref{Ansatz} into \eqref{PDE.uu} and discarding terms with prefactor $e^{\ii j(\kappa \cdot x -\omega t)/\epsi}$ if $|j|>3 $ yields the system
\begin{align}
\label{PDE.mfe}
\pt u_j + \frac{\ii}{\epsi}\L(j \omega, j \kappa)u_j + A(\partial)u_j
= \epsi \sum_{j_1+j_2+j_3=j} T(u_{j_1},u_{j_2},u_{j_3}), &
\\
\notag
\text{for } j\in\{1,3\}, \quad t\in (0,\tend/\epsi], \quad x\in\IR^d. &
\end{align}
The sum on the right-hand side is taken over the set 
\begin{align*}
 \Big\{J=(j_1,j_2,j_3) \in \J^3 : \#J:= j_1+j_2+j_3=j \Big\}.
\end{align*}
This set has only finitely many elements, namely 12 for $j=1$ and 10 for $j=3$.
Note that $ |J|_1 \geq |\#J| = |j_1 + j_2 + j_3 | = |j|$.
Since $u_{-j} = \overline{u_j}$ the PDEs for $u_{-1}$ and $u_{-3}$ 
are redundant but compatible.
\mynote{
The PDE \eqref{PDE.mfe} is compatible with the condition that $ u_{-j} = \overline{u_j} $, because 
\begin{align*}
 \overline{\frac{\ii}{\epsi}\L(j \omega, j \kappa)u_j} 
 &= 
 \frac{\ii}{\epsi}\L(-j \omega, -j \kappa)u_{-j},
 \\
 \sum_{j_1+j_2+j_3=j} \overline{T(u_{j_1},u_{j_2},u_{j_3})}
 &= \sum_{j_1+j_2+j_3=-j} T(u_{j_1},u_{j_2},u_{j_3}).
\end{align*}
F\"ur die zweite Gleichung braucht man, dass
$ \overline{T(a,b,c)}=T(\bar{a},\bar{b},\bar{c}) $ gilt. Daf\"ur braucht man keine zus\"atzliche Annahme, denn die Trilinearit\"at liefert
\begin{align*}
& T(a_1+\ii a_2,b_1+\ii b_2, c_1+\ii c_2) 
\\
&= 
T(a_1,b_1,c_1) 
+ \ii \Big(T(a_2,b_1,c_1) + T(a_1,b_2,c_1) + T(a_1,b_1,c_3) \Big)
\\
&\quad 
- \Big(T(a_2,b_2,c_1) + T(a_2,b_1,c_2) + T(a_1,b_2,c_2) \Big)
- \ii T(a_2,b_2,c_2),
\end{align*}
und damit gilt
\begin{align*}
& \overline{T(a_1+\ii a_2,b_1+\ii b_2, c_1+\ii c_2)} 
\\
&= 
T(a_1,b_1,c_1) 
-\ii \Big(T(a_2,b_1,c_1) + T(a_1,b_2,c_1) + T(a_1,b_1,c_3) \Big)
\\
&\quad 
- \Big(T(a_2,b_2,c_1) + T(a_2,b_1,c_2) + T(a_1,b_2,c_2) \Big)
+ \ii T(a_2,b_2,c_2)
\\
&= T(a_1-\ii a_2,b_1-\ii b_2, c_1-\ii c_2).
\end{align*}
}
The coupled system \eqref{PDE.mfe} is endowed with initial data
\begin{align}
\label{PDE.mfe.initial.data}
u_{\pm1}(0,\cdot) = u_{\pm1}^0 := p, \qquad u_{\pm3}(0,\cdot) = u_{\pm3}^0 = 0.
\end{align}

The main advantage of \eqref{PDE.mfe} over \eqref{PDE.uu} is that \eqref{PDE.mfe} does not oscillate in \emph{space}, because the initial data \eqref{PDE.mfe.initial.data} are smooth in contrast to 
\eqref{PDE.uu.b}. The price to pay is that the total number of unknowns in \eqref{PDE.mfe} is twice as large than in \eqref{PDE.uu}. Typical solutions of \eqref{PDE.mfe} still oscillate in \emph{time} due to the term $ \tfrac{\ii}{\epsi}\L(j \omega, j \kappa)u_j $, but it will turn out later that the situation is now more favourable; cf.~Remark~\ref{Remark.interpretation.bounds.u} below.

\mynote{An dieser Stelle sollte man wohl besser nicht mehr dazu sagen. Argumente wie ``Hauptteil im Kern'' sind gef\"ahrlich, weil wir ja sp\"ater nicht die Eigenzerlegung von $\L(\omega,\kappa)$ betrachten, sondern die von $ \ii \L(j \omega, j \kappa) + \epsi A(\partial)$ im Fourierraum. Dieser Operator hat keinen Kern, da r\"aumlich konstante Funktionen nicht in $L^1$ liegen.
Au\ss erdem oszillieren die $u_j$ ja tats\"achlich immer noch. Richtig gutartig wird das erst, wenn man die adiabatischen Variablen betrachtet bzw. die richtige Projektion verwendet. In Remark~\ref{Remark.interpretation.bounds.u} steht das etwas pr\"aziser.
}

Similar approximations have been considered in many other works.
A well-known approach in nonlinear geometric or diffractive optics is to look for an 
approximation of the form $ u(t,x) \approx \U(t,x,(\kappa \cdot x - \w t)/\epsi) $
for a profile $ \U=\U(t,x,\theta) $ which is periodic with respect to $\theta$; 
cf.~\cite{rauch-book:12,donnat-rauch:97b,donnat-joly-metivier-rauch:96,joly-metivier-rauch:98,colin-gallica-laurioux:05,lannes:98,alterman-rauch:03,joly-metivier-rauch:93}. In \cite{chartier-crouseilles-lemou-mehats-2015} this approach is used for the construction of uniformly accurate numerical methods for highly oscillatory problems.
The price to pay, however, is that introducing an additional variable increases the number of unknowns of the numerical discretization by a factor $ N_\theta$, where $N_\theta$ is the number of grid points in the $\theta$-direction. This is why we do not use a profile $ \U(t,x,\theta) $ explicitly.
However, the ansatz \eqref{Ansatz} can be interpreted as a truncated Fourier series
of $ \theta \mapsto \U(t,x,\theta) $, where all Fourier modes with index $ j\not\in\J$ are discarded.

Local well-posedness of the system \eqref{PDE.mfe} on long time intervals will be shown in Lemma~\ref{Lemma.wellposedness} below. 
Our main goal is to prove that under certain assumptions the approximation \eqref{Ansatz} has the accuracy
\begin{align*}
 \sup_{t\in[0,\tstar/\epsi]} \| \uu(t) - \uutilde(t) \|_{L^\infty(\IR^n)}
 &\leq C\epsi^2
\end{align*}
for $\epsi\in(0,1]$ with constants $\tstar$ and $C$ which do not depend on $\epsi$ (cf. Theorem~\ref{Theorem.error.bound}).
Comparing with \eqref{error.bound.nls} shows that $\uutilde$ is more accurate than
$\uu_{\text{\tiny NLS}}$.

\subsection{Evolution equations in Fourier space}

Let $\F f = \widehat{f} $ be the Fourier transform of $f$, i.e.
\begin{align*}
(\F f)(k) := (2\pi)^{-d/2}\intl_{\IR^d} f(x)e^{-\ii k\cdot x} dx
\end{align*}
with inverse Fourier transform 
\begin{align*}
f(x) = (2\pi)^{-d/2}\intl_{\IR^d} \widehat{f}(k)e^{\ii k\cdot x} dk.
\end{align*}
As in \cite{colin-lannes:09} we will work in the Wiener algebra  
\begin{align*}
W &= \{f \in \S'(\IR^d): \widehat{f}\in L^1(\IR^d)\}, &
\|f\|_W &= \|\widehat{f}\|_{L^1} = \intl_{\IR^d} |\widehat{f}(k)|_2 \; dk.
\end{align*}
$W(\IR^d)$ is a Banach algebra and continuously embedded in $L^\infty(\IR^d)$.
\mynote{Quelle: Gleichung (16) in Colin/Lannes.}
For $s\in\IN$ we set
\begin{align*}
W^s &= \{f \in W(\IR^d): \partial^\alpha f \in W(\IR^d) \text{ for all }
\alpha \in \IN_0^d, |\alpha|_1 \leq s\}, \\
\|f\|_{W^s} &= \sum_{|\alpha|_1 \leq s} \|\partial^\alpha f\|_W.
\end{align*}
Local well-posedness of \eqref{PDE.uu} in the Wiener algebra $W$ on time intervals 
$[0,\tend/\epsi]$ has been shown in \cite[Theorem 1]{colin-lannes:09}.

For estimates in the Wiener algebra it is convenient to consider the evolution equation \eqref{PDE.mfe}
in Fourier space. Applying the Fourier transform to the left-hand side of \eqref{PDE.mfe} gives
\begin{align*}
\F\Big(\pt u_j + \frac{\ii}{\epsi}\L(j\w,j\kappa)u_j + A(\partial)u_j \Big)(t,k)
&=
\pt \uhat_j(t,k) + \frac{\ii}{\epsi}\L_j(\epsi k)\uhat_j(t,k)
\end{align*}
with the shorthand notation 
\begin{align}
\label{Def.L.j}
\L_j(\theta) &:= \L(j\w, j\kappa + \theta) = -j\w I + A(j\kappa+\theta)-\ii E,
&& j\in\{1,3\}.
\end{align}
The Fourier transform of $T(u_{j_1},u_{j_2},u_{j_3})$ is given by
\begin{align}
\notag
\F\Big(T(u_{j_1},u_{j_2},u_{j_3})\Big)(k) &= (2\pi)^{-d} \intl_{\#K = k} 
T(\uhat_{j_1}(k^{(1)}),\uhat_{j_2}(k^{(2)}),\uhat_{j_3}(k^{(3)}))
\; dK
\\
\label{Def.T.Fourier}
&=: \T\left(\uhat_{j_1},\uhat_{j_2},\uhat_{j_3}\right)(k)
\end{align}
with $K=\big(k^{(1)},k^{(2)},k^{(3)}\big)\in \IR^d \times \IR^d \times \IR^d $, 
$\#K:=k^{(1)}+k^{(2)}+k^{(3)}\in \IR^d $, and with the notation
\begin{align*}
& \intl_{\#K = k} T(\uhat_{j_1}(k^{(1)}),\uhat_{j_2}(k^{(2)}),\uhat_{j_3}(k^{(3)}))
\; dK
\\
&=
\intl_{\IR^d} \intl_{\IR^d} 
T(\uhat_{j_1}(k^{(1)}),\uhat_{j_2}(k^{(2)}),\uhat_{j_3}(k-k^{(1)}-k^{(2)}))
\; dk^{(2)} \; dk^{(1)}.
\end{align*}
Hence, $u=(u_1,u_3)$ solves the system \eqref{PDE.mfe} if and only if 
$\uhat=(\uhat_1,\uhat_3)$ solves the system 
\begin{align}
\label{PDE.mfe.Fourier}
\pt \uhat_j(t,k) + \frac{\ii}{\epsi}\L_j(\epsi k)\uhat_j(t,k)
=
\epsi \sum_{\#J=j} \T\big(\uhat_{j_1},\uhat_{j_2},\uhat_{j_3}\big)(t,k), &
\\
\notag
j\in\{1,3\}, \quad t\in(0,\tend/\epsi], \quad k \in \IR^d &
\end{align}
with initial data
\begin{align}
\label{PDE.mfe.Fourier.initial.data}
\uhat_1(0,\cdot) = \phat, \qquad
\uhat_3(0,\cdot) = 0,
\end{align}
where $\phat$ is the Fourier transform of $p$ from \eqref{PDE.uu.b}.
The convention $u_{-j} = \overline{u_j}$ implies that $\uhat_{-j}(t,k) = \overline{\uhat_j(t,-k)}$.
For negative indices we set
\begin{align}
\label{Def.L.j.negative}
\L_{-j}(\theta)&= -\overline{\L_j(-\theta)},
&& j\in\{1,3\}
\end{align}
such that \eqref{PDE.mfe.Fourier} holds also for $j\in\{-1,-3\}$.

For later use, we note that \eqref{Def.T.Fourier} implies the inequality
\mynote{Stichwort ``algebra'' nennen?}
\begin{align}
\big\| \T(\widehat{f}_1,\widehat{f}_2,\widehat{f}_3) \big\|_{L^1}
\label{Bound.T.Fourier}
\leq 
C_\T \| \widehat{f}_1 \|_{L^1} \| \widehat{f}_2 \|_{L^1} \| \widehat{f}_3 \|_{L^1}
\end{align}
\mynote{
\begin{align}
\nonumber
 & \left\| \T(\widehat{f}_1,\widehat{f}_2,\widehat{f}_3) \right\|_{L^1}
 \\
\nonumber
&\leq
 (2\pi)^{-d} \intl_{\IR^d}  \intl_{\#K = k} 
\Big| T(\widehat{f}_1(k^{(1)}),\widehat{f}_2(k^{(2)}),\widehat{f}_3(k^{(3)})) \Big| \; dK \, dk
 \\
\nonumber
&\leq
 C_T (2\pi)^{-d} \intl_{\IR^d}  \intl_{\#K = k} 
| \widehat{f}_1(k^{(1)}) |_2 |\widehat{f}_2(k^{(2)})|_2 |\widehat{f}_3(k^{(3)}))|_2 \; dK \, dk
\\
\label{Bound.T.Fourier}
&\leq 
C_\T \| \widehat{f}_1 \|_{L^1} \| \widehat{f}_2 \|_{L^1} \| \widehat{f}_3 \|_{L^1}
\end{align}
}
for all $f_1,f_2,f_3 \in W(\IR^d)$, where $C_\T=C_T/(2\pi)^d$, and where $C_T$ is a constant such that
\begin{align*}
 | T(a,b,c) |_2 \leq C_T |a|_2 |b|_2 |c|_2 \qquad \text{for all } a,b,c\in\IC^n.
\end{align*}
\mynote{$T$ ist f\"ur Vektoren definiert, nicht f\"ur Funktionen. 
Da es auf dem $\IR^n$ nur endlich viele Basisvektoren gibt, gibt es ein $C_\star$ mit
$ |T(e_i,e_j,e_k)|_2 \leq C_\star $ f\"ur alle $i,j,k=1, \ldots,n$.
F\"ur beliebige Vektoren $a,b,c$ gilt dann 
\begin{align*}
 | T(a,b,c) |_2 
 &\leq \Big| \sum_{i} \sum_{j} \sum_{k} a_i b_j c_k  T(e_i,e_j,e_k) \Big|_2 
 \\
 &\leq C_\star |a|_1 |b|_1 |c|_1
 \leq C_\star n^{3/2} |a|_2 |b|_2 |c|_2
\end{align*}
wegen $|a|_1 \leq \sqrt{n} |a|_2$.
F\"ur $C_T = C_\star n^{3/2} $ gilt also die gew\"unschte Absch\"atzung.
}

\subsection{Local well-posedness on long time intervals}

For $ v=(v_1,v_3) \in W^s \times W^s$ we define the norm
\begin{align*}
 \| v \|_{W^s} =  2\| v_1 \|_{W^s} +  2\| v_3 \|_{W^s}.
\end{align*}
The factor $2$ is introduced in order to account for the terms
with negative indices which appear on the right-hand side of \eqref{PDE.mfe.Fourier}. If we set $\vhat_{-j}(k) = \overline{\vhat_j(-k)}$ for $j\in\{1,3\}$ as before, then
\begin{align*}
 \| v \|_{W^s} = \sum_{j\in\J} \| v_j \|_{W^s}.
\end{align*}
For $ v=(v_1,v_3) \in W \times W$ the inequalities
\begin{align}
\notag
\sum_{j\in\J}
 \Big\| \sum_{\#J=j} \T\big(\vhat_{j_1},\vhat_{j_2},\vhat_{j_3}\big) \Big\|_{L^1}
 &\leq 
 \sum_{J\in\J^3} \big\| \T\big(\vhat_{j_1},\vhat_{j_2},\vhat_{j_3}\big) \big\|_{L^1}
 \\
 \notag
  &\leq 
 C_\T \sum_{J\in\J^3} \| \vhat_{j_1} \|_{L^1} \| \vhat_{j_2} \|_{L^1} \| \vhat_{j_3} \|_{L^1}
 \\
\label{Bound.T.1}
 &\leq
 C_\T \| v \|_W^3 
\end{align}
follow from \eqref{Bound.T.Fourier}. If $ u=(u_1,u_3) $ is another element in 
$ W \times W$, then the trilinearity of $T$ and \eqref{Bound.T.Fourier} yield
\begin{align}
\notag
& \sum_{j\in\J} \sum_{\#J=j} \big\| 
\T\big(\uhat_{j_1},\uhat_{j_2},\uhat_{j_3}\big) 
- \T\big(\vhat_{j_1},\vhat_{j_2},\vhat_{j_3}\big) 
 \big\|_{L^1}
 \\
 \label{Bound.T.2}
 &\leq
 C_\T \Big(
  \| v \|_W^2 + 
  \| u \|_W \| v \|_W + 
  \| u \|_W^2 \Big)\| u-v \|_W.
\end{align}
After these preparations, local well-posedness of \eqref{PDE.mfe} can be shown.
The polarization of the initial data (Assumption~\ref{Ass:polarization}) is not required for the following result.

\begin{Lemma}[Local well-posedness]\label{Lemma.wellposedness}

\text{ }

\begin{enumerate}[label=(\roman*)]
\item
\label{Lemma.wellposedness.1}
If $ p \in W $, then there is a $\tend^\star>0$ such that for every $\epsi\in(0,1]$ the system \eqref{PDE.mfe} with initial data \eqref{PDE.mfe.initial.data}
has a unique mild solution 
\begin{align*}
u=(u_1,u_3), \qquad
u_j \in C([0,\tend^\star/\epsi),W).
\end{align*}

\item
\label{Lemma.wellposedness.2}
If $ p \in W^1 $ and $\tend<\tend^\star$, then the mild solution on 
$[0,\tend/\epsi]$ is a classical solution $u=(u_1,u_3)$ with
\begin{align*}
u_j \in C^1([0,\tend/\epsi],W) \cap C([0,\tend/\epsi],W^1).
\end{align*}

\item
\label{Lemma.wellposedness.3}
If $ p \in W^2 $ and $\tend<\tend^\star$, then 
\begin{align*}
u_j \in C^2([0,\tend/\epsi],W) \cap C^1([0,\tend/\epsi],W^1) \cap C([0,\tend/\epsi],W^2).
\end{align*}

\end{enumerate}
\end{Lemma}
By continuity, there are constants $C_{u,1}$ and $C_{u,2}$ such that
\begin{align}
\label{Lemma.wellposedness.bound.1}
\sup_{t\in[0,\tend/\epsi]}\| u_j(t) \|_{W^1} \leq C_{u,1},
\qquad j\in\{1,3\}
\end{align}
in case \ref{Lemma.wellposedness.2}, and 
\begin{align}
\label{Lemma.wellposedness.bound.2}
\sup_{t\in[0,\tend/\epsi]}\| u_j(t) \|_{W^2} \leq C_{u,2},
\qquad j\in\{1,3\}
\end{align}
in case \ref{Lemma.wellposedness.3}. In both cases the constant $C_{u,i}$ depends on $\tend$, $C_\T$ and on $\|p\|_{W^i}$, 
but not on $\epsi$.

The proof is based on classical arguments, but nevertheless we outline the main steps.

\proof
Choose a fixed $\epsi\in(0,1]$. The operator
\begin{align*}
\A\begin{pmatrix} v_1 \\ v_3 \end{pmatrix} &= 
\begin{pmatrix} \A_1 v_1 \\ \A_3 v_3 \end{pmatrix}, \qquad
\A_j = \frac{\ii}{\epsi}\L(j\omega,j\kappa) + A(\partial)
\end{align*}
with domain $ D(\A) = W^1 \times W^1 $ generates a strongly continuous group $ (e^{t\A})_{t\in\IR}$ on $ W \times W $.
For $j\in\{1,3\}$ and every $t\in\IR$ the group operator $  e^{t\A_j} $ is an isometry, because 
\begin{align*}
\| e^{t\A_j}v_j \|_W &= \big\| \F\big(e^{t\A_j}v_j\big) \big\|_{L^1}
= \intl_{\IR^d} \big| e^{\ii t \L_j(\epsi k)/\epsi}\vhat_j(k) \big|_2 \; dk 
\\
&= \intl_{\IR^d} \big| \vhat_j(k) \big|_2 \; dk = \| v_j \|_W
\end{align*}
for all $v_j \in W$ since $ \L_j(\epsi k) $ is Hermitian.
The system \eqref{PDE.mfe} can be reformulated as 
\begin{align}
 \label{PDE.mfe.Cauchy.problem}
\pt u_j + \A_j u_j = \epsi \sum_{\#J=j} T(u_{j_1},u_{j_2},u_{j_3}) 
\qquad 
\text{for } j\in\{1,3\} .
\end{align}
For a number $ \tau>0 $ to be determined below, consider the space
\begin{align*}
\X = C([0,\tau/\epsi],W) \times C([0,\tau/\epsi],W)
 \end{align*}
with norm
\begin{align*}
\| v \|_\X = \sup_{t\in[0,\tau/\epsi]} \|v(t)\|_W
 = \sup_{t\in[0,\tau/\epsi]} \sum_{j\in\J} \|v_j(t)\|_W 
\end{align*}
and the mapping 
\begin{align*}
 & \Phi \colon \X \longrightarrow \X, \qquad 
 \Phi(v) = v^{new}=\begin{pmatrix} v_1^{new} \\ v_3^{new} \end{pmatrix},
\\
&
v_j^{new}(t) = e^{-t\A_j}u_j^0 + \epsi \sum_{\#J=j} \intl_0^t e^{(s-t)\A_j} 
T(v_{j_1},v_{j_2},v_{j_3})(s) \; ds.
\end{align*}
With \eqref{Def.T.Fourier} we obtain
\begin{align*}
\| v_j^{new}(t) \|_W 
&\leq 
\| u_j^0 \|_W + \epsi \sum_{\#J=j} \intl_0^t \| T(v_{j_1},v_{j_2},v_{j_3})(s) \|_W \; ds
\\
&=
\| u_j^0 \|_W + \epsi \sum_{\#J=j} \intl_0^t \| \T(\vhat_{j_1},\vhat_{j_2},\vhat_{j_3})(s) \|_{L^1} \; ds
\end{align*}
and with \eqref{Bound.T.1} it follows that
\begin{align*}
\| \Phi(v) \|_\X 
&= \sup_{t\in[0,\tau/\epsi]} \sum_{j\in\J} \|v_j^{new}(t)\|_W
\\
&\leq 
2\| p \|_W + 
\epsi \sup_{t\in[0,\tau/\epsi]} \sum_{j\in\J}
\sum_{\#J=j} \intl_0^t \| \T(\vhat_{j_1},\vhat_{j_2},\vhat_{j_3})(s) \|_{L^1} \; ds
\\
&\leq 
2\| p \|_W + 
C_\T \epsi  \sup_{t\in[0,\tau/\epsi]} \intl_0^t 
\| v(s) \|_W^3 \; ds
\\
&\leq 
2\| p \|_W + 
C_\T \tau  \sup_{s\in[0,\tau/\epsi]} \| v(s) \|_W^3 
\\
&= 
2\| p \|_W + 
C_\T \tau \|v\|_\X^3.
\end{align*}
Now choose $\rho>0$ and set $r=1+\rho$. 
For every $ p\in W $ with $ \| p \|_W\leq \tfrac{\rho}{2} $, $\Phi$ maps the ball
\begin{align*}
 B(r)=\{v \in \X \colon \| v \|_\X \leq r\}
\end{align*}
onto itself under the condition that $ \tau \leq 1/(C_\T r^3) $.
If $ v,w \in B(r) $,
then it follows from \eqref{Bound.T.2} that
\begin{align*}
\| \Phi(v) - \Phi(w) \|_\X
&= 
\epsi \sup_{t\in[0,\tau/\epsi]}
\sum_{j\in\J} \sum_{\#J=j} \intl_0^t 
\| T(v_{j_1},v_{j_2},v_{j_3})(s) - T(w_{j_1},w_{j_2},w_{j_3})(s) \|_W
\; ds
\\
&= 
3 C_\T r^2 \epsi \sup_{t\in[0,\tau/\epsi]} \intl_0^t \| v(s)-w(s) \|_W \; ds
\\
&\leq
3C_\T r^2 \tau \| v - w \|_\X.
\end{align*} 
If we choose, e.g., $ \tau = \min\{1/(C_\T r^3), 1/(6C_\T r^2)\}$,
then $ \Phi: B(r) \rightarrow B(r) $ is a contraction, and by
Banach's fixed point theorem, there is a unique fixed point $u \in B(r)$ of $ \Phi$. By construction, $u$ is a mild solution of 
\eqref{PDE.mfe} with initial data \eqref{PDE.mfe.initial.data}.
With standard arguments, this solution can be extended to a maximal time interval
$ [0,\tau^+(\epsi)/\epsi) $, and one can show that
\mynote{Das folgt daraus, dass $\tau^+(\epsi) \geq \tau $ unabh\"angig von $\epsi$ nach unten beschr\"ankt ist und man daher auch $\tend^\star \geq \tau $ hat, oder?}
\begin{align*}
\tend^\star := \inf_{\epsi\in(0,1]}\tau^+(\epsi)>0.
\end{align*}
This proves part \ref{Lemma.wellposedness.1}.

The nonlinearity $T:W\times W \times W \rightarrow W $ is continuously differentiable
due to \eqref{Bound.T.Fourier}. 
\mynote{Die Frechet-Ableitung von $T(u,v,w)$ ist
\begin{align*}
 h \mapsto T(u,v,h)+T(u,h,w)+T(h,v,w).
\end{align*}
Diese Abbildung ist stetig in $u,v,w$ wegen \eqref{Bound.T.Fourier}.}
If $ p \in W^1 $, then $ (u_1^0,u_3^0)=(p,0) \in D(\A) $,
and applying Theorem 1.5 in \cite[Chapter 6]{pazy:83} yields that for every $\tend<\tend^\star$
the mild solution is in fact a classical solution on $[0,\tend/\epsi]$, which proves part \ref{Lemma.wellposedness.2}.

To show part \ref{Lemma.wellposedness.3} we set $u'=(u_1',u_3')$ with $u'_j = \pt u_j$ 
and formally differentiate both sides of \eqref{PDE.mfe.Cauchy.problem} to obtain
\begin{subequations}
 \label{PDE.mfe.Cauchy.problem.derivative}
\begin{align}
\pt u_j' + \A_j u_j' = \epsi \sum_{\#J=j} 
\Big(T(u_{j_1}',u_{j_2},u_{j_3}) + T(u_{j_1},u_{j_2}',u_{j_3}) +T(u_{j_1},u_{j_2},u_{j_3}')\Big)
\end{align}
with initial data
\begin{align}
 \label{PDE.mfe.Cauchy.problem.derivative.initial.data}
u_j'(0) &= -\A_j u_j^0 + \epsi \sum_{\#J=j} T(u_{j_1}^0,u_{j_2}^0,u_{j_3}^0).
\end{align}
\end{subequations}
Let $ u=(u_1,u_3)$ be the classical solution constructed in part \ref{Lemma.wellposedness.2} and consider the linear problem
\begin{subequations}
\label{PDE.mfe.Cauchy.problem.derivative.B}
\begin{align}
\pt u_j' + \A_j u_j' &= \epsi \B_j(t,u'),
\\
\notag
\B_j(t,u') &= \sum_{\#J=j} 
\Big(T(u_{j_1}'(t),u_{j_2}(t),u_{j_3}(t)) + T(u_{j_1}(t),u_{j_2}'(t),u_{j_3}(t)) 
\\
&\hs{20}
+ T(u_{j_1}(t),u_{j_2}(t),u_{j_3}'(t))\Big)
\end{align}
\end{subequations}
with initial data \eqref{PDE.mfe.Cauchy.problem.derivative.initial.data}.
Since $ u_j \in C^1([0,\tend/\epsi],W) \cap C([0,\tend/\epsi],W^1) $ the mapping
\begin{align*}
(t,u') \mapsto \B_j(t,u'), \qquad \B_j \colon [0,\tend/\epsi] \times (W\times W) \rightarrow W 
\end{align*}
is continuously differentiable, and if $ p \in W^2 $, then $ (u_1'(0),u_3'(0)) \in D(\A) $ due to \eqref{Bound.T.Fourier}.
Hence, the mild solution $ u'=\pt u $ of \eqref{PDE.mfe.Cauchy.problem.derivative.B}, \eqref{PDE.mfe.Cauchy.problem.derivative.initial.data} is in fact a classical solution according to Theorem 1.5 in \cite[Chapter 6]{pazy:83}, which proves part \ref{Lemma.wellposedness.3}.
\qed

\subsection{Transformation to smoother variables}\label{Subsec.Transformation}

In order to analyze the accuracy of the approximation \eqref{Ansatz}, 
the estimate \eqref{Lemma.wellposedness.bound.1} has to be refined. 
In Section~\ref{Sec:Refined.bounds} we will show that if the system \eqref{PDE.mfe} is considered with initial data \eqref{PDE.mfe.initial.data}, then $ u_3 $ stays small on long time intervals, i.e.
\begin{align*}
\sup_{t\in[0,\tend/\epsi]}\| u_3(t) \|_{W^1} &\leq C\epsi.
\end{align*}
A similar estimate will be shown for a certain ``part'' of $u_1$ to be specified later; cf.~\eqref{Bound.uhat.1}.
In order to formulate and prove these refined estimates, it is very useful to consider a transformation of $\uhat_j$ which is introduced now.

For $j\in\{\pm1, \pm3\}$ and every $\theta\in\IR^d$ the Hermitian matrix 
$ \L_j(\theta) = \L(j\w, j\kappa + \theta)$ defined in \eqref{Def.L.j} has an eigendecomposition
\begin{align*}
 \L_j(\theta) &= \Psi_j(\theta)\Lambda_j(\theta)\Psi_j^*(\theta) 
\end{align*}
with a unitary matrix\footnote{In contrast to the traditional notation we do not denote the unitary matrix by $U_j(\theta)$ in order to avoid confusion with 
$ \uu $ from \eqref{PDE.uu}, $u_j$ from \eqref{PDE.mfe}, $U_{\text{\tiny NLS}}$,
or $U_{\text{\tiny SVEA}}$.} 
$ \Psi_j(\theta)\in\IC^{n \times n} $ and a real diagonal matrix
$ \Lambda_j(\theta)\in\IR^{n \times n} $ containing the eigenvalues of  $ \L_j(\theta) $.
It follows from \eqref{Def.L.j.negative} that
$\Psi_{-j}(\theta)= -\overline{\Psi_j(-\theta)}$ and 
$\Lambda_{-j}(\theta)= -\overline{\Lambda_j(-\theta)}= -\Lambda_j(-\theta)$.
Let $\psi_{j\ell}(\theta) \in \IC^n $ be the $\ell$-th column of 
 $ \Psi_j(\theta) $, and let $\lambda_{j\ell}(\theta) \in \IR $ be the $\ell$-th eigenvalue, i.e. 
\begin{align}
\label{eigenbasis}
 \L_j(\theta)\psi_{j\ell}(\theta) &= 
 \lambda_{j\ell}(\theta)\psi_{j\ell}(\theta),
 &
 \psi_{jm}(\theta) \cdot \psi_{j\ell}(\theta) &= 
 \begin{cases} 1 & \text{if } m=\ell, \\ 0 & \text{else} \end{cases} 
\end{align}
for $m,\ell=1, \ldots, n$. 
The matrix $ \L_1(0)=\L(\w,\kappa) $ has a one-dimensional kernel according to 
Assumption~\ref{Ass:polarization}\ref{Ass:polarization.item.i},
and the enumeration of the eigenvalues is chosen in such a way that $ \lambda_{11}(0) = 0$.
Hence, the kernel of $\L_1(0)$ is spanned by $\psi_{11}(0)$, which is important in the context of Assumption~\ref{Ass:polarization}\ref{Ass:polarization.item.ii}.

For every $j\in\{1,3\},\epsi>0,t\geq 0,k\in\IR^d$ we define 
\begin{align}
  \label{Def.TT}
 \TT_{j,\epsi}(t,k)&= \exp\big(\tfrac{\ii t}{\epsi}\Lambda_{j}(\epsi k)\big)
 \Psi_{j}^*(\epsi k)
 = \Psi_{j}^*(\epsi k) \exp\big(\tfrac{\ii t}{\epsi}\L_{j}(\epsi k)\big)
\end{align}
and consider the new variables
\begin{align}
 \label{Def.z}
 z_{j}(t,k) &= \TT_{j,\epsi}(t,k)\uhat_j(t,k),
 \qquad j\in\{1,3\},
\end{align}
where $ \uhat_1(t,k), \uhat_3(t,k) $ is the solution of \eqref{PDE.mfe.Fourier}.
The matrix \eqref{Def.TT} is unitary, and hence 
\begin{align}
\label{Norm.uhat.z}
| \uhat_j(t,k) |_2 = | z_j(t,k) |_2
\quad \text{and} \quad 
 \| \uhat_j(t) \|_{L^1} = \| z_j(t) \|_{L^1}.
\end{align}
For negative indices we set
\begin{align*}
\TT_{-j,\epsi}(t,k) &:= \overline{\TT_{j,\epsi}(t,-k)}, &
z_{-j}(t,k) &:= \overline{z_{j}(t,-k)}.
\end{align*}
Taking the time derivative of \eqref{Def.z}, substituting 
\eqref{PDE.mfe.Fourier} and using
\begin{align*}
\pt \TT_{j,\epsi}(t,k)
= \tfrac{\ii}{\epsi}\Lambda_{j}(\epsi k) \TT_{j,\epsi}(t,k)
= \frac{\ii}{\epsi}\TT_{j,\epsi}(t,k)\L_j(\epsi k)
\end{align*}
yields
\begin{align}
\label{zdot.in.terms.of.uhat}
\pt z_{j}(t) 
&=
\epsi \sum_{\#J=j} \FF(t,\uhat,J)(t)
\end{align}
with $\uhat=(\uhat_1,\uhat_3)$ and
\begin{align}
\label{Def.F.uhat}
\FF(t,\uhat,J) = \TT_{j,\epsi}(t) \T\big(\uhat_{j_1},\uhat_{j_2},\uhat_{j_3}\big)(t), \qquad j=\#J.
\end{align}
A closed system of evolution equations for $ z_1 $ and $z_3$ could be obtained by 
expressing the right-hand side of \eqref{Def.F.uhat} by $ z_j$ via  
the inverse transform $ \uhat_j(t,k) = \TT_{j,\epsi}^*(t,k)z_{j}(t,k) $. 
Since this leads to rather complicated formulas, we will avoid this whenever possible.

Comparing \eqref{PDE.mfe.Fourier} with \eqref{zdot.in.terms.of.uhat}
shows that the dominating linear term
\begin{align*}
\frac{\ii}{\epsi}\L_j(\epsi k)\uhat_j(t,k)
\end{align*}
in \eqref{PDE.mfe.Fourier} is removed by the transformation.
In the \emph{linear} case (where $\T(\cdot,\cdot,\cdot)=0$) it follows from \eqref{zdot.in.terms.of.uhat}
and \eqref{Def.F.uhat} that $ \pt z_j(t) = 0$, and hence that $ z_j(t)=z_j(0) $ is constant in time.
The exact solution of \eqref{PDE.mfe.Fourier} is then simply 
\begin{align*}
\uhat_j(t,k) 
&= \TT_{j,\epsi}^*(t,k)\TT_{j,\epsi}(0,k)\uhat_j(0,k) 
= \TT_{j,\epsi}^*(t,k)z_j(0,k) \qquad
 \text{for } \T(\cdot,\cdot,\cdot)=0. 
\end{align*}
This favourable property does not cure the oscillatory behaviour completely
in the general (nonlinear) case, but 
the entries of $z_{j}$ oscillate with a much smaller amplitude than
the entries of $\uhat_j$. The reason is that the right-hand side of \eqref{zdot.in.terms.of.uhat}
is formally $\Ord{\epsi}$ instead of $\Ord{1/\epsi}$ in \eqref{PDE.mfe.Fourier}. This is our main motivation for considering transformed variables in the proofs of our main results.

Initial data for $z_1$ and $z_3$ are obtained from \eqref{PDE.mfe.Fourier.initial.data}
and \eqref{Def.z}, namely
\begin{align}
\label{PDE.z.initial.data}
z_{j}(0,k) &= \TT_{j,\epsi}(0,k) \uhat_j(0,k)
  =
  \begin{cases}
  \Psi_1^*(\epsi k) \phat(k) & \text{if } j=1, \\ 0 & \text{if } j=3.
   \end{cases}
\end{align}
Since $\text{ker}\big(\L(\w,\kappa)\big)=\text{ker}\big(\L_1(0)\big)
= \text{span}\{\psi_{11}(0)\}$ it follows from Assumption~\ref{Ass:polarization} 
that
$\psi_{1\ell}^*(0) \widehat{p}(k)=0$ for all $k$ and all $\ell \not=1$, but this is in general not true for $\psi_{1\ell}^*(\epsi k) \widehat{p}(k)$.
For initial data $ p \in W^1 $, however, it has been shown in \cite[proof of Lemma 3, page 718]{colin-lannes:09}
that
\begin{align*}
\| \psi_{1\ell}(\epsi \cdot)\psi_{1\ell}^*(\epsi \cdot) \widehat{p} \|_{L^1} \leq C \epsi
\| \nabla p \|_{W} 
\end{align*}
for every $\ell \not=1$, which yields
\begin{align}
\label{Ass:polarization.itemz}
 \| z_{1\ell}(0) \|_{L^1} \leq C \epsi \|  \nabla p \|_{W} 
 \qquad \text{for all } \ell \not=1.
\end{align}
The special role of the first entry of $z_1$ can be expressed by means of the projection
\begin{align}
\label{Def.P1}
P: \IC^n \rightarrow \IC^n, \quad
(w_1, \ldots, w_n)^T \mapsto (w_1, 0, \ldots, 0)^T.
\end{align}
With $\Pperp=(I-P)$, Assumption~\ref{Ass:polarization} implies that
$ \| \Pperp z_1(0) \|_{L^1} = \Ord{\epsi}$, because
\begin{align}
\notag
 \| \Pperp z_1(0) \|_{L^1} 
&=
 \intl_{\IR^d} | \Pperp z_1(0,k) |_2 \; dk
 \leq
 \sum_{\ell=2}^n \| z_{1\ell}(0) \|_{L^1}
 \\
 \label{Bound.projection.initial.value}
 &\leq C (n-1) \epsi \|  \widehat{\nabla p} \|_{L^1} 
 \end{align}
\mynote{
\begin{align}
\notag
 \| \Pperp z_1(0) \|_{L^1} 
 &=
 \intl_{\IR^d} | \Pperp z_1(0,k) |_2 \; dk
 \\
\notag
 &\leq
 \intl_{\IR^d} | \Pperp z_1(0,k) |_1 \; dk
 \\
\notag
 &=
 \sum_{\ell=2}^n \intl_{\IR^d} | z_{1\ell}(0,k) | \; dk
 \\
\notag
 &=
 \sum_{\ell=2}^n \| z_{1\ell}(0) \|_{L^1}
 \\
 \label{Bound.projection.initial.value}
 &\leq C (n-1) \epsi \|  \widehat{\nabla p} \|_{L^1} 
 \end{align}
 } 
due to \eqref{Ass:polarization.itemz}.
For later use, we also define the projections
\begin{align}
 \label{Def.P2}
 \widehat{w} \mapsto \P_\epsi \widehat{w}, \qquad
\P_\epsi(k)\widehat{w}(k) = \psi_{11}(\epsi k)\psi_{11}^*(\epsi k)\widehat{w}(k)
\end{align}
and $\P_\epsi^\perp = I-\P_\epsi.$
$\P_\epsi$ projects a vector-valued function $ \widehat{w} : \IR^d \rightarrow \IC^n$ pointwise into the first eigenspace of $\L_1(\epsi k)$.
For $ \uhat_1(t,k) = \TT_{j,\epsi}^*(t,k) z_1(t,k) $ it follows that
\begin{align}
\label{Relation.between.projections}
\P_\epsi(k)\uhat_1(t,k) = \TT_{j,\epsi}^*(t,k) Pz_1(t,k),
\end{align}
which means, in particular, that 
\begin{align*}
\P_\epsi(0)\uhat_1(t,0) \in \text{ker}\big(\L_1(0)\big) = \text{span}\{\psi_{11}(0)\}.
\end{align*}

\section{Refined bounds for the coefficient functions}\label{Sec:Refined.bounds}

\subsection{Setting and goal}\label{Subsec.setting}

Let $u=(u_1,u_3)$ with 
$u_j \in C^1([0,\tend/\epsi],W) \cap C([0,\tend/\epsi],W^1)$ be a classical
solution of \eqref{PDE.mfe} with initial data \eqref{PDE.mfe.initial.data}
for some $ p \in W^1 $. 
Our next goal is to prove that under certain assumptions 
there is a $\tstar\in(0,\tend]$ independent of $\epsi$ and a constant $C$ such that 
\begin{subequations}
\label{Bounds.uhat}
\begin{align}
\label{Bound.uhat.1}
\sup_{t\in[0,\tstar/\epsi]}\| \P_\epsi^\perp \uhat_1(t) \|_{L^1} 
&\leq C \epsi,
\\
\label{Bound.uhat.3}
\sup_{t\in[0,\tstar/\epsi]}\| \uhat_3(t) \|_{L^1} 
&\leq C \epsi
\end{align}
\end{subequations}
for all $\epsi\in(0,1].$
For $z_{j} = \TT_{j,\epsi}\uhat_j$ these bounds are equivalent to
\begin{subequations}
\label{Bounds.z}
\begin{align}
\label{Bound.z.1}
\sup_{t\in[0,\tstar/\epsi]}\| \Pperp z_1(t) \|_{L^1} 
&\leq C \epsi,
\\
\label{Bound.z.3}
\sup_{t\in[0,\tstar/\epsi]}\| z_3(t) \|_{L^1} 
&\leq C \epsi
\end{align}
\end{subequations}
for all $\epsi\in(0,1]$ due to \eqref{Norm.uhat.z} and \eqref{Relation.between.projections}.
\mynote{Wegen der Projektion kann man diese Aussagen nur f\"ur $\uhat_j$ formulieren, nicht f\"ur $u_j$.}

\begin{Remark}
 \label{Remark.interpretation.bounds.u}
The bounds \eqref{Bounds.uhat} and \eqref{Bounds.z} are not only crucial for proving error bounds for the approximation \eqref{Ansatz} (cf.~Theorem~\ref{Theorem.error.bound} below),
but also interesting from a numerical point of view.
In fact, \eqref{Bounds.uhat} means that
$ \uhat_1(t) = \P_\epsi \uhat_1(t) + \Ord{\epsi} $
and $ \uhat_3(t) = \Ord{\epsi} $, such that the ``main part'' of the solution of \eqref{PDE.mfe.Fourier} is $\P_\epsi \uhat_1(t)$. But this part is \emph{essentially non-oscillatory} according to Lemma~\ref{Lemma.bound.uhat1dot} below.
This can be exploited in the construction of efficient numerical methods for \eqref{PDE.mfe.Fourier}.
We are currently working on the analysis of such methods.

\mynote{Die zweite Ableitung $\P_\epsi \pt^2 \uhat_1(t,k)$ ist immer noch $\Ord{1}$, weil man zwar beim Durchdifferenzieren ein $1/\epsi$ erh\"alt, was aber durch das $\epsi$ vor der Nichtlinearit\"at kompensiert wird. Allerdings ist unklar, ob man eine Funktion $y(t)$ mit $\dot{y}(t)=\Ord{1}$ auch dann als nichtoszillatorisch bezeichnen kann, wenn man alles auf einem Intervall der L\"ange $\Ord{1/\epsi}$ betrachtet. Dieses Fass m\"ochte ich an dieser Stelle nicht aufmachen -- ich denke das ist ok so. Im Fall von $ \P_\epsi \uhat_1(t,k) $ kann man an dieser Stelle ja ausnutzen, dass man die $\Ord{1}$-Dynamik durch das mitbewegte Koordinatensystem eliminieren kann. H\"ohere Ableitungen von $ \P_\epsi \uhat_1(t,k) $ werden nat\"urlich Oszillationen haben -- deshalb habe ich oben ``essentially non-oscillatory'' geschrieben.
}

\end{Remark}

In order to prove \eqref{Bounds.z} we define the scaled norm
\begin{align}
 \label{Def.scaled.norm}
  \il y \il_\epsi = 2\| Py_1 \|_{L^1} + \frac{2}{\epsi} \| \Pperp y_1 \|_{L^1}
  + \frac{2}{\epsi} \| y_3 \|_{L^1}
\end{align}
for all $ y = (y_1,y_3) $ with $ y_j \in L^1(\IR^d,\IC^n)$.
As before, we set $ y_{-j} = \overline{y_j} $.
Since
\begin{align}
\label{Replace.scaled.norm.01}
 \| y_1 \|_{L^1} \leq \| Py_1 \|_{L^1} + \| \Pperp y_1 \|_{L^1}
 \leq \| Py_1 \|_{L^1} + \epsi^{-1} \| \Pperp y_1 \|_{L^1}
\end{align}
holds for all $ \epsi\in(0,1]$, it follows that
\begin{align}
 \label{Replace.scaled.norm.02}
  \sum_{j\in\J} \epsi^{(1-|j|)/2} \| y_j \|_{L^1} 
  &\leq \il y \il_\epsi,
\end{align}
which will be used frequently. Our goal is to prove that there is a constant $C$ such that
\begin{align*}
\sup_{t\in[0,\tstar/\epsi]} \il z(t) \il_\epsi \leq C,
\end{align*}
for all $\epsi\in(0,1]$, because this implies \eqref{Bounds.z} and hence also \eqref{Bounds.uhat}.

\begin{Proposition}\label{Proposition01}
Let $u=(u_1,u_3)$ be the classical
solution of \eqref{PDE.mfe} with initial data \eqref{PDE.mfe.initial.data}
for some $ p \in W^1 $.
Let $z_1$ and $z_3$ be the transformed variables defined in \eqref{Def.z}.
For every sufficiently large $ r>0$ there is a $\tstar\in(0,\tend]$ such that 
under Assumptions~\ref{Ass:polarization} and \ref{Ass:L.properties} 
\begin{align*}
\sup_{t\in[0,\tstar/\epsi]} \il z(t) \il_\epsi \leq r \qquad
\text{for all } \epsi\in(0,1].
\end{align*}
The constant $\tstar$ depends on $\tend$, $r$, $C_{u,1}$, $C_\T$, on the inverse of the nonzero eigenvalues of $\Lambda_1(0)$ and $\Lambda_3(0)$, and on the Lipschitz constant in Assumption \ref{Ass:L.properties}\ref{Ass:L.properties.item.ii}, but not on $\epsi$.
\end{Proposition}

\begin{Remark}
 ``Sufficiently large'' means that $ r $ must be larger than the constant $C_\bullet$ which occurs in the proof. This condition is required to ensure that $ \tstar$ defined in \eqref{Prop01.tstar.formula} is positive.
\end{Remark}

\mynote{Die Konstante $C_\bullet$ h\"angt unter anderem auch von $\tend$ ab. Das sieht widerspr\"uchlich aus, weil man $z(t)$ ja nur auf dem \emph{kleineren} Zeitintervall $ [0,\tstar/\epsi] \subset [0,\tend/\epsi] $ betrachtet. Trotzdem macht das Sinn, denn man kann $\tstar$ nur dann durch die Formel \eqref{Prop01.tstar.formula} definieren, wen die Konstanten $ \widehat{C}$ und $C_\bullet$ unabh\"angig von $\tstar$ sind. Deshalb kann man es wohl nicht vermeiden, $\tstar$ manchmal durch $\tend$ zu ersetzen.
Das wollen wir aber nicht erl\"autern -- das stiftet nur Verwirrung.
}

\mynote{Reminder: $|v|_2 \leq |v|_1 \leq \sqrt{n} |v|_2 $. Das sollte man dort erw\"ahnen, wo es zum ersten mal verwendet wird.}

The following lemmas will be used in the proof of Proposition~\ref{Proposition01}.


\begin{Lemma}\label{Lemma.useful.bounds}
Let $\vhat=(\vhat_1,\vhat_3)$ with $\vhat_j\in L^1$ and $ \vhat_{-j}(k)=\overline{\vhat_j(-k)}$ for $j=1,3$. 
For every $J=(j_1,j_2,j_3)\in\J^3$ the inequality
 \begin{align*}
\big\| \FF(t,\vhat,J) \big\|_{L^1} 
&\leq 
C_\T \prod_{i=1}^3 \| \TT_{j_i,\epsi}(t)\vhat_{j_i} \|_{L^1},
 \end{align*}
holds for all $t\geq 0$.
\end{Lemma}

\proof
The definition \eqref{Def.F.uhat}, the inequality \eqref{Bound.T.Fourier} and the fact that $\TT_{j,\epsi}$ is unitary 
imply that
\begin{align*}
\big\| \FF(t,\vhat,J) \big\|_{L^1} 
\leq
\left\| \T(\vhat_{j_1},\vhat_{j_2},\vhat_{j_3}) \right\|_{L^1}
\leq 
C_\T \| \vhat_{j_1} \|_{L^1} \| \vhat_{j_2} \|_{L^1} \| \vhat_{j_3} \|_{L^1}.
\end{align*}
Now the assertion follows from $\| \vhat_{j_i} \|_{L^1} = \| \TT_{j_i,\epsi}(t)\vhat_{j_i} \|_{L^1}.$
\qed


The PDE system \eqref{PDE.mfe.Fourier} suggests that formally $ \pt \uhat_j(t) = \Ord{1/\epsi} $.
The following lemma shows, however, that $ \pt \P_\epsi \uhat_1(t) $ can be bounded independently of $\epsi$ 
on long time intervals.

 \begin{Lemma}\label{Lemma.bound.uhat1dot}
Under the assumptions of Lemma~\ref{Lemma.wellposedness} \ref{Lemma.wellposedness.2}
there is a constant $C$ such that
 \begin{align*}
\sup_{t\in[0,\tend/\epsi]} 
\| \pt \P_\epsi \uhat_1(t) \|_{L^1} \leq C.
\end{align*}
$C$ depends on the constant $C_{u,1}$ from \eqref{Lemma.wellposedness.bound.1} and thus also on $\tend$, but not on $\epsi$.
\end{Lemma}
This lemma corresponds to Lemma 2 in \cite{colin-lannes:09}. The proof is based on the observation that
\begin{align}
\notag
 \frac{\ii}{\epsi}\P_\epsi \L_1(\epsi k)\uhat_1(t,k)
 &=
 \frac{\ii}{\epsi}\lambda_{11}(\epsi k)\P_\epsi \uhat_1(t,k)
 \end{align}
because of \eqref{Def.P2} and \eqref{eigenbasis}. The Lipschitz continuity of the eigenvalues 
(cf.~Assumption~\ref{Ass:L.properties}\ref{Ass:L.properties.item.ii})) and the fact that $\lambda_{11}(0)=0$
yield
\begin{align*}
| \lambda_{11}(\epsi k) | = | \lambda_{11}(\epsi k) - \lambda_{11}(0) | \leq C\epsi |k|_1,
\end{align*}
and thus
\begin{align*}
\Big| \frac{\ii}{\epsi}\P_\epsi \L_1(\epsi k)\uhat_1(t,k) \Big|_2
\leq 
C |k|_1 |\P_\epsi \uhat_1(t,k)|_2.
\end{align*}
Together with \eqref{PDE.mfe.Fourier} this shows that $\P_\epsi \pt \uhat_1(t,k)=\pt \P_\epsi \uhat_1(t,k)$ is uniformly bounded.

\subsection{Proof of Proposition \ref{Proposition01}}
\label{Subsec.proof.prop01}

According to \eqref{zdot.in.terms.of.uhat} and \eqref{Def.scaled.norm} we have
\begin{align}
 \notag
\il z(t) \il_\epsi 
&\leq 
 \il z(0) \il_\epsi + \il \intl_0^t \pt z(s) \; ds \il_\epsi 
 \\
 \notag
 & \leq
 \il z(0) \il_\epsi
 + 2 \sum_{\# J=1} \left(
 \epsi \Big\| \intl_0^t P \FF(s,\uhat,J) \; ds \Big\|_{L^1}
 +
 \Big\| \intl_0^t \Pperp \FF(s,\uhat,J) \; ds \Big\|_{L^1}
 \right)
 \\
 \label{Prop01.bound10}
 & \quad + 
 2 \sum_{\# J=3} \Big\| \intl_0^t \FF(s,\uhat,J) \; ds \Big\|_{L^1}
\end{align}
with $\FF$ defined in \eqref{Def.F.uhat}. The first term 
\begin{align*}
 \il z(0) \il_\epsi 
 &= 2\| Pz_1(0) \|_{L^1} + \frac{2}{\epsi} \| \Pperp z_1(0) \|_{L^1} + \frac{2}{\epsi} \| z_3(0) \|_{L^1}
\end{align*}
is uniformly bounded, because $ z_3(0)=0 $ according to \eqref{PDE.z.initial.data}
and $\| \Pperp z_1(0) \|_{L^1} \leq C\epsi \|p\|_{W^1}$ due to \eqref{Bound.projection.initial.value}.
Now let
\begin{align}
\label{Def.aj}
a_j(t) &=
\begin{cases}
\| Pz_1(t) \|_{L^1} + \epsi^{-1} \| \Pperp z_1(t) \|_{L^1} 
& \text{if } j=\pm 1, \\
\epsi^{-1} \| z_3(t) \|_{L^1} & \text{if } j=\pm 3,
\end{cases}
\end{align}
which according to \eqref{Def.scaled.norm} means that
\begin{align}
\label{Scaled.norm.sum.aj}
\sum_{j\in\J} a_j(t) = 2a_1(t) + 2a_3(t) = \il z(t) \il_\epsi.
\end{align}
We will prove that there are constants $C_\star$ and $\Chat$ such that
\begin{align}
\label{Target.term.1}
\epsi \Big\| \intl_0^t P \FF(s,\uhat,J) \; ds \Big\|_{L^1}
+
\Big\| \intl_0^t \Pperp \FF(s,\uhat,J) \; ds \Big\|_{L^1}
&
\\
\notag
\leq
C_\star +  \frac{\Chat \epsi}{2} \intl_0^t \prod_{i=1}^3 a_{j_i}(s) \; ds
&
\end{align}
holds for all $t\in[0,\tend/\epsi]$ and for every $J=(j_1,j_2,j_3)\in\J^3$ with $\#J=1$, and that
\begin{align}
\label{Target.term.3}
\Big\| \intl_0^t \FF(s,\uhat,J) \; ds \Big\|_{L^1}
\leq
C_\star +  \frac{\Chat \epsi}{2} \intl_0^t \prod_{i=1}^3 a_{j_i}(s) \; ds
\end{align}
holds for all $t\in[0,\tend/\epsi]$ and for every $J\in\J^3$ with $\#J=3$.
If this is true, then substituting into \eqref{Prop01.bound10} yields
\begin{align*}
 \il z(t) \il_\epsi 
  &\leq 
  C_\bullet
  + \Chat \epsi \sum_{j\in\J} \sum_{\# J=j} 
 \intl_0^t \prod_{i=1}^3 a_{j_i}(s) \; ds
 \\
  &\leq
  C_\bullet + \Chat \epsi 
 \intl_0^t \bigg(\sum_{j\in\J} a_j(s)\bigg)^3 \; ds
 \\
 &=
C_\bullet + \Chat \epsi \intl_0^t \il z(s) \il_\epsi^3 \; ds
\end{align*}
by \eqref{Scaled.norm.sum.aj} with a constant $C_\bullet$ which depends on $\|p\|_{W^1}$, $C_\star$ and the (finite) number of multi-indices $J$ with $\#J=1$ and $\#J=3$, respectively. Now let 
\begin{align*}
\sigma_\epsi(t)=\sup_{s\in[0,t]}\il z(s) \il_\epsi
\end{align*}
and observe that 
\begin{align}
\label{Prop01.before.condition.tstar}
\sigma_\epsi(t) 
&\leq
  C_\bullet + \Chat \epsi \intl_0^t \il z(s) \il_\epsi^3 \; ds
  \leq
  C_\bullet + \Chat \epsi t \sigma_\epsi^3(t).
\end{align}
Let $r>C_\bullet$. 
If we can choose $\tstar$ in such a way that 
\begin{align*}
C_\bullet + \Chat \tstar \sigma_\epsi^3(\tstar) \leq r,
\end{align*}
then \eqref{Prop01.before.condition.tstar} and the fact that $\sigma_\epsi$ is monotonically increasing
implies that 
$ \sigma_\epsi(t) \leq r $ for all $t\in[0,\tstar/\epsi]$.
Hence, we choose
\begin{align}
\label{Prop01.tstar.formula}
  \tstar = \frac{r-C_\bullet}{\Chat r^3}.
\end{align}
Of course, the choice \eqref{Prop01.tstar.formula} is in most cases way too pessimistic.
What is important is that $\tstar$ depends on $C_\bullet$, $r$, and $\Chat$, but not on $\epsi$.

Now the inequalities \eqref{Target.term.1} and \eqref{Target.term.3} have to be shown. This is the main part of the proof.

\subsubsection{Proof of \eqref{Target.term.3}}
\label{Subsubsec.proof.prop01.term.3}

Let $J\in\J$ with $\# J=3$. Two cases have to be treated separately.

\paragraph{Case 1: \boldmath $|J|_1>3$.}
In this case we have $|J|_1 \geq 5$ because $|J|_1$ is odd.
As an example, the reader may consider $J=(1,3,-1)$.
Lemma~\ref{Lemma.useful.bounds} 
and the fact that 
\begin{align*}
\epsi^{(|J|_1-3)/2}\prod_{i=1}^3 \epsi^{(1-|j_i|)/2} = \epsi^0 = 1 
\end{align*}
yields
\begin{align*}
  \Big\| \intl_0^t \FF(s,\uhat,J) \; ds\Big\|_{L^1}
  &\leq
  C_\T  
  \intl_0^t 
  \prod_{i=1}^3 \| z_{j_i}(s) \|_{L^1}  
  \; ds  
  \\
  &=
  C_\T \epsi^{(|J|_1-3)/2} 
  \intl_0^t 
  \prod_{i=1}^3 \Big( \epsi^{(1-|j_i|)/2} \| z_{j_i}(s) \|_{L^1} \Big)
  \; ds  
  \\
  &\leq
  C_\T \epsi \intl_0^t \prod_{i=1}^3 a_{j_i}(s) \; ds,
\end{align*}
because $ \epsi^{(|J|_1-3)/2} \leq \epsi $ and 
$\epsi^{(1-|j_i|)/2} \| z_{j_i}(s) \|_{L^1} \leq a_{j_i}(s)$
by definition \eqref{Def.aj}. 
This yields an estimate of the form \eqref{Target.term.3} with $C_\star=0$ and
$\Chat=2 C_\T$.

\paragraph{Case~2: \boldmath $|J|_1=3 $.}
This case appears only if $J=(1,1,1)$.
\mynote{Im Fall $J=-(1,1,1)$ gilt $ \#J=-3 \not= 3$.}
In this situation, the simple argument from Case 1 is not enough to prove the desired bound, because now $ \epsi^{(|J|_1-3)/2} = 1$.
One power of $\epsi$ has to be gained from the oscillatory behavior of 
$\FF(s,\uhat,J)$.
First, the nonlinearity in
\begin{align*}
\Big\| \intl_0^t \FF(s,\uhat,J) \; ds\Big\|_{L^1}
&=
\Big\| \intl_0^t 
  \TT_{3,\epsi}(s) \T(\uhat_1,\uhat_1,\uhat_1)(s)
\; ds\Big\|_{L^1}
\end{align*}
is split into eight parts:
\begin{align}
\label{T.uhat1.decomp}
 \T\big(\uhat_1,\uhat_1,\uhat_1\big)
 &= \T\big(\P_\epsi\uhat_1, \P_\epsi\uhat_1, \P_\epsi\uhat_1\big) 
 + \T\big(\P_\epsi \uhat_1, \P_\epsi \uhat_1, \P_\epsi^\perp \uhat_1\big)
 \\
 \notag
 &\quad + \T\big(\P_\epsi\uhat_1, \P_\epsi^\perp\uhat_1, \P_\epsi\uhat_1\big) 
 + \T\big(\P_\epsi \uhat_1, \P_\epsi^\perp\uhat_1, \P_\epsi^\perp\uhat_1\big)
 \\
 \notag
 &\quad + \quad \ldots \quad + \T\big(\P_\epsi^\perp \uhat_1, \P_\epsi^\perp\uhat_1, \P_\epsi^\perp\uhat_1\big).
\end{align}
All terms where $\P_\epsi^\perp \uhat_1$ appears in at least one of the three arguments are easy to treat.
With \eqref{Bound.T.Fourier} and \eqref{Norm.uhat.z} we obtain for example
\begin{align}
\notag
& \Big\| \intl_0^t 
\TT_{3,\epsi}(s) \T\big(\P_\epsi^\perp \uhat_1,\P_\epsi\uhat_1,\P_\epsi\uhat_1\big)(s)
\; ds\Big\|_{L^1}
\leq
\intl_0^t
\big\| \T\big(\P_\epsi^\perp \uhat_1,\P_\epsi\uhat_1,\P_\epsi\uhat_1\big)(s) \big\|_{L^1} \; ds
\\
\notag
&\leq
C_\T \, \epsi \intl_0^t \Big(
\tfrac{1}{\epsi}\| \P_\epsi^\perp \uhat_1(s) \|_{L^1} 
\cdot \| \P_\epsi\uhat_1(s) \|_{L^1} 
\cdot \| \P_\epsi\uhat_1(s) \|_{L^1}
\Big)
\; ds
\\
\notag
&=
C_\T \, \epsi \intl_0^t \Big(
\tfrac{1}{\epsi}\| \Pperp z_1(s) \|_{L^1} 
\cdot \| Pz_1(s) \|_{L^1} 
\cdot \| Pz_1(s) \|_{L^1}
\Big)
\; ds
\\
\label{Proof.Prop01.bound02.case.2.easy.part}
&\leq
C_\T \, \epsi 
\intl_0^t \prod_{i=1}^3 a_{j_i}(s) \; ds,
\end{align}
because \eqref{Def.aj} implies that 
$ \| z_1(s) \|_{L^1} \leq a_1(s)$. All other parts of \eqref{T.uhat1.decomp} with the exception of 
$\T\big(\P_\epsi\uhat_1, \P_\epsi\uhat_1, \P_\epsi\uhat_1\big)$ can be treated in the same way, i.e. for each of these terms we obtain an estimate of the type
\eqref{Target.term.3} with $C_\star=0$ and $\Chat=2 C_\T$.

The main difficulty in Case~2 is to prove that
\begin{align}
\label{Proof.Prop01.bound02.case.2.difficult.term}
& \Big\| \intl_0^t 
\TT_{3,\epsi}(s)\T\big(\P_\epsi \uhat_1,\P_\epsi \uhat_1,\P_\epsi \uhat_1\big)(s)
\; ds\Big\|_{L^1} \leq C
\end{align}
uniformly in $\epsi$ in spite of the integration over a possibly long time interval.
\eqref{Proof.Prop01.bound02.case.2.difficult.term} corresponds to a bound of the form \eqref{Target.term.3},
but this time with $\Chat=0$.
By \eqref{Def.TT} we have
\begin{align*}
\TT_{3,\epsi}(s,k) &= 
\exp\big(\tfrac{\ii s}{\epsi}\Lambda_{3}(0)\big) 
\exp\big(\tfrac{\ii s}{\epsi}\Delta_{3}(\epsi k)\big) 
\Psi_{3}^*(\epsi k)
\end{align*}
with $ \Delta_{3}(\epsi k) = \Lambda_{3}(\epsi k) - \Lambda_{3}(0) $.
Hence, \eqref{Proof.Prop01.bound02.case.2.difficult.term} can be expressed as
\begin{align*}
& \Big\| \intl_0^t 
\exp\big(\tfrac{\ii s}{\epsi}\Lambda_{3}(0)\big) f_\epsi(s)
\; ds\Big\|_{L^1},
\\
f_\epsi(t,k) &=
\exp\big(\tfrac{\ii t}{\epsi}\Delta_{3}(\epsi k)\big) 
\Psi_{3}^*(\epsi k)\T\big(\P_\epsi \uhat_1,\P_\epsi \uhat_1,\P_\epsi \uhat_1\big)(t,k).
\end{align*}
Since the matrix
$ \L_3(0) = \L(3\w, 3\kappa) $ is invertible by Assumption~\ref{Ass:polarization}\ref{Ass:polarization.item.iii}, we can integrate by parts and obtain
\begin{align*}
& \Big\| \intl_0^t 
\exp\big(\tfrac{\ii s}{\epsi}\Lambda_{3}(0)\big) f_\epsi(s)
\; ds\Big\|_{L^1}
\\
&\leq
\Big\| \Big[
\tfrac{\epsi}{\ii}\Lambda_{3}^{-1}(0)
\exp\big(\tfrac{\ii s}{\epsi}\Lambda_{3}(0)\big) f_\epsi(s)
\Big]_0^t
\Big\|_{L^1}
+
\Big\| 
\tfrac{\epsi}{\ii}\Lambda_{3}^{-1}(0)
\intl_0^t 
\exp\big(\tfrac{\ii s}{\epsi}\Lambda_{3}(0)\big) \pt f_\epsi(s)
\; ds\Big\|_{L^1}
\\
&\leq
C\epsi \Big(\| f_\epsi(0) \|_{L^1} + \| f_\epsi(t) \|_{L^1}\Big)
+
C\epsi \intl_0^t 
\Big\| \pt f_\epsi(s) \; ds\Big\|_{L^1}.
\end{align*}
With $|\TT_{3,\epsi}(s,k)|_2=1$, \eqref{Bound.T.Fourier}, and \eqref{Lemma.wellposedness.bound.1} it follows that
\begin{align*}
 \| f_\epsi(t) \|_{L^1}
 &=
  \big\| \T\big(\P_\epsi \uhat_1,\P_\epsi \uhat_1,\P_\epsi \uhat_1\big)(t) \big\|_{L^1}
\leq C_\T \| \P_\epsi \uhat_1(t) \|_{L^1}^3 
\leq C
\end{align*}
with a constant which depends on $C_\T$ and the constant $ C_{u,1}$ from \eqref{Lemma.wellposedness.bound.1}.
$ \Lambda_3$ is globally Lipschitz continuous by Assumption~\ref{Ass:L.properties}\ref{Ass:L.properties.item.ii}, i.e.
\begin{align*}
 |\tfrac{\ii}{\epsi}\Delta_{3}(\epsi k)|_2 = \tfrac{1}{\epsi}|\Lambda_{3}(\epsi k) - \Lambda_{3}(0)|_2 \leq C|k|_1
\end{align*}
with a constant $C$ which does not depend on $\epsi$ and $k$.
This yields for $t\in[0,\tend/\epsi]$ that
\begin{align*} 
\epsi \intl_0^t 
\Big\| \pt f_\epsi(s) \; ds\Big\|_{L^1}
& \leq
\tend \sup_{s\in[0,\tend/\epsi]}
\Big\| \pt f_\epsi(s) \Big\|_{L^1}
\\
&\leq
C \tend \sup_{s\in[0,\tend/\epsi]}
\intl_{\IR^d} |k|_1 \cdot | \T\big(\P_\epsi \uhat_1,\P_\epsi \uhat_1,\P_\epsi \uhat_1\big)(s,k) |_2 \; dk
\\
&\quad +
C \tend \sup_{s\in[0,\tend/\epsi]}
\intl_{\IR^d} | \pt \T\big(\P_\epsi \uhat_1,\P_\epsi \uhat_1,\P_\epsi \uhat_1\big)(s,k) |_2 \; dk
\\
&\leq
C C_\T \tend \Big(C_{u,1}^3 + C_{u,1}^2 \sup_{s\in[0,\tend/\epsi]} \| \pt \P_\epsi \uhat_1(s) \|_{L^1}
\Big)
\end{align*}
by \eqref{Bound.T.Fourier} and \eqref{Lemma.wellposedness.bound.1}.
Since $\sup_{s\in[0,\tend/\epsi]} \| \pt \P_\epsi \uhat_1(s) \|_{L^1}$ is uniformly bounded according to 
Lemma~\ref{Lemma.bound.uhat1dot}, it follows that \eqref{Proof.Prop01.bound02.case.2.difficult.term}
is uniformly bounded. All in all this proves the inequality \eqref{Target.term.3} in Case~2.

\subsubsection{Proof of \eqref{Target.term.1}}
\label{Subsubsec.proof.prop01.term.1}

Let $J\in\J^3$ with $\#J=1$. With Lemma~\ref{Lemma.useful.bounds} and the fact that $ \epsi^{(1-|j_i|)/2} \geq 1 $ for $|j_i|\geq 1$ it follows that
\begin{align*}
\notag
\epsi \Big\| \intl_0^t P\FF(s,\uhat,J) \; ds \Big\|_{L^1} 
&\leq 
\epsi \intl_0^t \big\| \FF(s,\uhat,J) \big\|_{L^1} \; ds 
\\
\notag
&\leq 
C_\T \, \epsi   \intl_0^t 
\prod_{i=1}^3 \| z_{j_i}(s) \|_{L^1} \; ds
\\
\notag
&\leq 
C_\T \, \epsi   \intl_0^t 
\prod_{i=1}^3 \Big(\epsi^{(1-|j_i|)/2} \| z_{j_i}(s) \|_{L^1} \Big) \; ds
\\
&\leq 
C_\T \, \epsi  \intl_0^t \prod_{i=1}^3 a_{j_i}(s) \; ds,
\end{align*}
because by definition $\epsi^{(1-|j_i|)/2} \| z_{j_i}(s) \|_{L^1} \leq a_{j_i}(s).$
This is a bound of type \eqref{Target.term.3} with $C_\star=0$.
Finally, we have to show a bound for the term
\begin{align*}
\Big\| \intl_0^t \Pperp \FF(s,\uhat,J) \; ds  \Big\|_{L^1} 
\end{align*}
for all $J\in\J^3$ with $\#J=1$. 
For $|J|_1 \geq 5$ we can simply proceed as in Case~1 in \ref{Subsubsec.proof.prop01.term.3},
but if $|J|_1 = 3$, then one power of $\epsi$ has to be gained.
The only multi-indices $J\in\J^3$ with $\#J=1$ and $|J|_1=3$ are 
$(1,1,-1)$, $(1,-1,1)$, and $(-1,1,1)$. 
We consider $J=(1,1,-1)$, because the other two permutations can be treated in the same way.
 
We can use \eqref{Def.F.uhat} to obtain
\begin{align*}
\Big\| \intl_0^t \Pperp \FF(s,\uhat,J) \; ds\Big\|_{L^1}
&=
\Big\| \intl_0^t 
  \Pperp \TT_{1,\epsi}(s) \T(\uhat_1,\uhat_1,\uhat_{-1})(s)
\; ds\Big\|_{L^1}.
\end{align*}
With $ \uhat_1 = \P_\epsi\uhat_1 + \P_\epsi^\perp \uhat_1 $ and
$ \uhat_{-1} = \overline{\P_\epsi}\uhat_{-1} + \overline{\P_\epsi^\perp} \uhat_{-1} $
the nonlinearity is split into eight parts
\begin{align*}
 \T\big(\uhat_1,\uhat_1,\uhat_{-1}\big)
 &= 
 \T\big(\P_\epsi\uhat_1,\P_\epsi\uhat_1,\overline{\P_\epsi}\uhat_{-1}\big)
 + 
 \T\big(\P_\epsi\uhat_1,\P_\epsi\uhat_1,\overline{\P_\epsi^\perp}\uhat_{-1}\big)
 \\
 &\quad + 
 \T\big(\P_\epsi\uhat_1,\P_\epsi^\perp\uhat_1,\overline{\P_\epsi}\uhat_{-1}\big)
 + 
 \T\big(\P_\epsi\uhat_1,\P_\epsi^\perp\uhat_1,\overline{\P_\epsi^\perp}\uhat_{-1}\big)
 \\
 & \quad
 + \quad \ldots \quad + 
 \T\big(\P_\epsi^\perp\uhat_1,\P_\epsi^\perp\uhat_1,\overline{\P_\epsi^\perp}\uhat_{-1}\big)
\end{align*}
similar to \eqref{T.uhat1.decomp}.
All terms containing $ \P_\epsi^\perp\uhat_1 $ or $ \overline{\P_\epsi^\perp}\uhat_{-1} $ 
can again be estimated in a straightforward way. The only remaining term
\begin{align}
\label{Proof.Prop01.bound01.difficult.term}
\Big\| \intl_0^t 
  \Pperp \TT_{1,\epsi}(s) 
  \T\big(\P_\epsi\uhat_1,\P_\epsi\uhat_1,\overline{\P_\epsi}\uhat_{-1}\big)(s)
\; ds\Big\|_{L^1}
\end{align}
has to be treated in a similar way as \eqref{Proof.Prop01.bound02.case.2.difficult.term}.
We let $ \Delta_{1}(\epsi k) = \Lambda_{1}(\epsi k) - \Lambda_{1}(0) $ and obtain from
\eqref{Def.TT} that
\begin{align*}
\Pperp \TT_{1,\epsi}(s,k) &= 
\exp\big(\tfrac{\ii s}{\epsi}\Lambda_{1}(0)\big) 
\Pperp \exp\big(\tfrac{\ii s}{\epsi}\Delta_{1}(\epsi k)\big) 
\Psi_{1}^*(\epsi k),
\end{align*}
because $\Pperp$ commutes with every diagonal matrix.
Hence, \eqref{Proof.Prop01.bound01.difficult.term} can be expressed as
\begin{align*}
& \Big\| \intl_0^t 
\exp\big(\tfrac{\ii s}{\epsi}\Lambda_{1}(0)\big) \Pperp f_\epsi(s)
\; ds\Big\|_{L^1},
\\
f_\epsi(t,k) &=
\exp\big(\tfrac{\ii t}{\epsi}\Delta_{1}(\epsi k)\big) \Psi_{1}^*(\epsi k)
\T\big(\P_\epsi\uhat_1,\P_\epsi\uhat_1,\overline{\P_\epsi}\uhat_{-1}\big)(t,k).
\end{align*}
In order to gain the missing factor $\epsi$ we want to integrate by parts. The diagonal matrix
$\Lambda_{1}(0)=\diag\big(\lambda_{11}(0), \ldots, \lambda_{1n}(0)\big)$ is not invertible, because 
$ \lambda_{11}(0)=0 $ (cf.~Section~\ref{Subsec.Transformation}), but this is compensated by the projection $ \Pperp $ which sets the first entry of a vector to zero. Hence, we can simply replace $ \lambda_{11}(0)$ by 1 (or any other nonzero number)
and consider $\widetilde{\Lambda}_{1}(0)=\diag\big(1,\lambda_{12}(0), \ldots, \lambda_{1n}(0)\big)$ instead of 
$\Lambda_{1}(0)$. The modified matrix $\widetilde{\Lambda}_{1}(0)$ is invertible because $ \lambda_{1\ell}(0)\not=0 $ for $\ell>1$ by Assumption~\ref{Ass:polarization}\ref{Ass:polarization.item.i}.
Integrating by parts yields
\begin{align*}
& \Big\| \intl_0^t 
\exp\big(\tfrac{\ii s}{\epsi}\widetilde{\Lambda}_{1}(0)\big) \Pperp f_\epsi(s)
\; ds\Big\|_{L^1}
\\
&\leq
\Big\| \Big[
\tfrac{\epsi}{\ii}\widetilde{\Lambda}_{1}^{-1}(0)
\exp\big(\tfrac{\ii s}{\epsi}\widetilde{\Lambda}_{1}(0)\big) \Pperp f_\epsi(s)
\Big]_0^t
\Big\|_{L^1}
+
\Big\| 
\tfrac{\epsi}{\ii}\widetilde{\Lambda}_{1}^{-1}(0)
\intl_0^t 
\exp\big(\tfrac{\ii s}{\epsi}\widetilde{\Lambda}_{1}(0)\big) \Pperp \pt f_\epsi(s)
\; ds\Big\|_{L^1}
\\
&\leq
C\epsi \Big(\| f_\epsi(0) \|_{L^1} + \| f_\epsi(t) \|_{L^1}\Big)
+
C\epsi \intl_0^t 
\Big\| \pt f_\epsi(s) \; ds\Big\|_{L^1}.
\end{align*}
As in Case 2 in Section~\ref{Subsubsec.proof.prop01.term.3} 
one can show that these terms are bounded by a constant which does not depend on $\epsi$.
This shows that \eqref{Proof.Prop01.bound01.difficult.term} is uniformly bounded.
Together with the other considerations, this proves the inequality \eqref{Target.term.1} and completes the proof of Proposition~\ref{Proposition01}.
\qed

\subsection{Extension to a stronger norm}

For $\mu\in\{1, \ldots d\}$ let $ D_\mu $ denote the Fourier multiplicator $(D_\mu \widehat{w})(k) = \ii k_\mu \widehat{w}(k).$
If $u=(u_1,u_3)$ is the classical solution of \eqref{PDE.mfe} and \eqref{PDE.mfe.initial.data}, then by
definition $ D_\mu \uhat_j  = \TT_{j,\epsi}^* D_\mu z_{j} $ is the Fourier transform of $\partial_\mu u_j$. 
(Note that the scalar multiplicator $D_\mu$ commutes with the matrix $\TT_{j,\epsi}^*$.)

In the proof of Theorem~\ref{Theorem.error.bound} we will also need the following version of Proposition~\ref{Proposition01} 
where $ z_j(t) $ is replaced by $ D_\mu z_j(t) $.

\begin{Proposition}\label{Proposition02} 
Let $u=(u_1,u_3)$ be the classical
solution of \eqref{PDE.mfe} with initial data \eqref{PDE.mfe.initial.data}
for some $ p \in W^2 $.
Let $z_1$ and $z_3$ be the transformed variables defined in \eqref{Def.z}, and let 
$\mu\in\{1, \ldots, d\}$.
Under the assumptions of Proposition~\ref{Proposition01} there is a constant $C$ such that
\begin{align*}
\sup_{t\in[0,\tstar/\epsi]} \il D_\mu z(t) \il_\epsi \leq C \qquad
\text{for all } \epsi\in(0,1]
\end{align*}
with $\tstar$ from Proposition~\ref{Proposition01}.
$C$ depends on $ \| p \|_{W^2} $, $C_{u,2}$ from \eqref{Lemma.wellposedness.bound.2}, and on $r$ from Proposition~\ref{Proposition01}, but not on $\epsi$.
\end{Proposition}

Proposition~\ref{Proposition02} implies the bounds
\begin{subequations}
\label{Bounds.Dmu.uhat}
\begin{align}
\label{Bound.Dmu.uhat.1}
\sup_{t\in[0,\tstar/\epsi]}\| D_\mu \P_\epsi^\perp \uhat_1(t) \|_{L^1} 
&\leq C \epsi,
\\
\label{Bound.Dmu.uhat.3}
\sup_{t\in[0,\tstar/\epsi]}\| D_\mu \uhat_3(t) \|_{L^1} 
&\leq C \epsi
\end{align}
\end{subequations}
for all $\epsi\in(0,1].$
For $z_{j} = \TT_{j,\epsi}\uhat_j$ these bounds are equivalent to
\begin{subequations}
\label{Bounds.Dmu.z}
\begin{align}
\label{Bound.Dmu.z.1}
\sup_{t\in[0,\tstar/\epsi]}\| D_\mu \Pperp z_1(t) \|_{L^1} 
&\leq C \epsi,
\\
\label{Bound.Dmu.z.3}
\sup_{t\in[0,\tstar/\epsi]}\| D_\mu z_3(t) \|_{L^1} 
&\leq C \epsi
\end{align}
\end{subequations}
for all $\epsi\in(0,1]$ due to \eqref{Norm.uhat.z} and \eqref{Relation.between.projections}.

\bigskip


\proofof{Proposition \ref{Proposition02}}
Let $\mu\in\{1, \ldots, d\}$ be fixed. Applying the operator $D_\mu$ to both sides of \eqref{zdot.in.terms.of.uhat}
gives
\begin{align}
\label{PDE.Dz}
\pt D_\mu z_{j}(t) 
&=
\epsi \sum_{\#J=j} D_\mu \FF(t,\uhat,J)(t), \qquad j\in\{1,3\}.
\end{align}
Integrating \eqref{PDE.Dz} from $0$ to $t$ and applying the scaled norm \eqref{Def.scaled.norm} leads to
\begin{align*}
 \il D_\mu z(t) \il_\epsi 
&\leq 
 \il D_\mu z(0) \il_\epsi + \il \intl_0^t \pt D_\mu z(s) \; ds \il_\epsi 
 \\
 & \leq
 \il D_\mu z(0) \il_\epsi 
 + 2 \sum_{\# J=1} f_1(t,\epsi,\uhat,J)
 + 2 \sum_{\# J=3} f_3(t,\epsi,\uhat,J)
 \intertext{with}
 f_1(t,\epsi,\uhat,J)
 &=
 \epsi \Big\| \intl_0^t P D_\mu \FF(s,\uhat,J) \; ds \Big\|_{L^1}
 +
 \Big\| \intl_0^t \Pperp D_\mu \FF(s,\uhat,J) \; ds \Big\|_{L^1},
 \\
 f_3(t,\epsi,\uhat,J)
 &=
 \Big\| \intl_0^t D_\mu \FF(s,\uhat,J) \; ds \Big\|_{L^1}.
\end{align*}
According to \eqref{PDE.z.initial.data} and \eqref{Bound.projection.initial.value} the term $ \il D_\mu z(0) \il_\epsi $ is uniformly bounded with a constant which depends on $\|p\|_{W^2}$.
\mynote{Anstelle von \eqref{Bound.projection.initial.value} verwenden wir 
$ \| \Pperp D_\mu z_1(0) \|_{L^1} \leq C\epsi \| p \|_{W^2} $.}
Our goal is to prove that there are constants $C_1$ and $C_2$ such that the inequality
\begin{align}
f_j(t,\epsi,\uhat,J)
 &\leq 
 C_1 + C_2 \, \epsi \intl_0^t \il D_\mu z(s) \il_\epsi  \; ds
 \label{Gronwall.bound}
\end{align}
holds for $j\in\{1,3\}$ and for all $ J\in\J^3 $ with $\#J=j. $
If \eqref{Gronwall.bound} is true, then it follows that
\begin{align*}
 \il D_\mu z(t) \il_\epsi 
\leq 
c_1 + c_2 \, \epsi \intl_0^t \il D_\mu z(s) \il_\epsi  \; ds
\end{align*}
(with other constants), and Gronwall's lemma yields
\begin{align*}
 \sup_{t\in[0,\tstar/\epsi]} \il D_\mu z(t) \il_\epsi \leq c_1 e^{c_2 \tstar},
\end{align*}
which proves the assertion.

To prove \eqref{Gronwall.bound} we analyze the term $D_\mu \FF(t,\uhat,J)$, which appears in $f_j(t,\epsi,\uhat,J)$.
By \eqref{Def.F.uhat} we have
\begin{align}
\label{Def.DF}
D_\mu \FF(t,\uhat,J) &= \TT_{j,\epsi}(t) D_\mu \T\big(\uhat_{j_1},\uhat_{j_2},\uhat_{j_3}\big)(t), \qquad j=\#J,
\end{align}
and it follows from \eqref{Def.T.Fourier} that
\begin{align}
\label{Product.rule.DT}
D_\mu \T\big(\uhat_{j_1},\uhat_{j_2},\uhat_{j_3}\big)
&=
\T\big(D_\mu\uhat_{j_1},\uhat_{j_2},\uhat_{j_3}\big)
+ \T\big(\uhat_{j_1},D_\mu\uhat_{j_2},\uhat_{j_3}\big) + \T\big(\uhat_{j_1},\uhat_{j_2},D_\mu\uhat_{j_3}\big),
\end{align}
which corresponds to the product rule. If we consider $ \uhat=(\uhat_1,\uhat_3)$ as given, then
\eqref{Product.rule.DT} is \emph{linear} with respect to $ D_\mu \uhat_j = \TT_{j,\epsi}^* D_\mu z_j $.
This is the reason why Gronwall's lemma can be used to prove boundedness of $ \il D_\mu z(t) \il_\epsi $,
but not to show boundedness of $ \il z(t) \il_\epsi $ in the proof of Proposition~\ref{Proposition01}.

Since $ \TT_{j,\epsi} $ is unitary it follows from \eqref{Def.DF} and \eqref{Product.rule.DT} that
\begin{align*}
& \big\| D_\mu \FF(t,\uhat,J) \big\|_{L^1}
=
\big\| D_\mu \T\big(\uhat_{j_1},\uhat_{j_2},\uhat_{j_3}\big)(t) \big\|_{L^1}
\\
&\leq
\| \T\big(D_\mu\uhat_{j_1},\uhat_{j_2},\uhat_{j_3}\big)(t) \|_{L^1}
+ \| \T\big(\uhat_{j_1},D_\mu\uhat_{j_2},\uhat_{j_3}\big)(t) \|_{L^1} 
+ \| \T\big(\uhat_{j_1},\uhat_{j_2},D_\mu\uhat_{j_3}\big)(t) \|_{L^1}.
\end{align*}
With \eqref{Bound.T.Fourier} we obtain for all $t\in[0,\tstar/\epsi]$ that
\begin{align}
\notag
\big\| D_\mu \FF(t,\uhat,J) \big\|_{L^1}
&\leq
C_\T
\Big(
\| D_\mu \uhat_{j_1}(t) \|_{L^1} \| \uhat_{j_2}(t) \|_{L^1} \| \uhat_{j_3}(t) \|_{L^1} 
\\
\notag
& \hs{15}
+ \| \uhat_{j_1}(t) \|_{L^1} \| D_\mu \uhat_{j_2}(t) \|_{L^1} \| \uhat_{j_3}(t) \|_{L^1} 
\\
\notag
& \hs{15}
+ \| \uhat_{j_1}(t) \|_{L^1} \| \uhat_{j_2}(t) \|_{L^1} \| D_\mu \uhat_{j_3}(t) \|_{L^1} 
\Big)
\\
\notag
&\leq 
C_\T C \Big(
\epsi^{(|j_2|+|j_3|-2)/2} \| D_\mu z_{j_1}(t) \|_{L^1}
+ \epsi^{(|j_1|+|j_3|-2)/2} \| D_\mu z_{j_2}(t) \|_{L^1}
\\
\label{Bound.DetaFj}
& \hs{15} 
+ \epsi^{(|j_1|+|j_2|-2)/2} \| D_\mu z_{j_3}(t) \|_{L^1}
\Big),
 \end{align}
because we already know that
\begin{align}
\label{Bounds.z.no.projection}
\sup_{t\in[0,\tstar/\epsi]}\| \uhat_j(t) \|_{L^1} = \sup_{t\in[0,\tstar/\epsi]}\| z_j(t) \|_{L^1} &\leq C \epsi^{(|j|-1)/2}
\end{align}
by Proposition \ref{Proposition01} and \eqref{Replace.scaled.norm.02}.

Let $\#J=j=3$. If $|J|_1 \geq 5$ we can proceed in a straightforward way: we estimate
\begin{align*}
\big\|f_3(t,\epsi,\uhat,J) \big\|_{L^1}
\leq 
\intl_0^t \big\| D_\mu \FF(s,\uhat,J) \big\|_{L^1} \; ds,
\end{align*}
substitute \eqref{Bound.DetaFj} and use that for $ |J|_1\geq 5$ there is at least
one index $j_i=3$. It is enough to consider $J=(3,1,1)$, because
$ J=(1,3,1)$, $J=(1,1,3)$ and all permutations of $J=(3,3,-3)$ can be treated in the same way. 
Combining \eqref{Bound.DetaFj} with \eqref{Bounds.z.no.projection} yields for $J=(j_1,j_2,j_3)=(3,1,1)$
\begin{align*}
\big\|f_3(t,\epsi,z,J) \big\|_{L^1}
&\leq 
C \intl_0^t 
\Big(
\| D_\mu z_3(s) \|_{L^1}
+ 2 \epsi \| D_\mu z_1(s) \|_{L^1}
\Big)
\; ds
\\
&\leq
C \epsi \intl_0^t \il D_\mu z(s) \il_\epsi  \; ds,
\end{align*}
which corresponds to \eqref{Gronwall.bound} with $C_1=0$.
Now let $|J|_1=\#J=j=3$, i.e. let $ J=(1,1,1)$. Then it follows from 
\eqref{Def.DF} and \eqref{Product.rule.DT} that
\begin{align}
\notag
\big\|f_3(t,\epsi,\uhat,J) \big\|_{L^1}
&=
\Big\| \intl_0^t 
  \TT_{3,\epsi}(s) \T(D_\mu \uhat_1,\uhat_1,\uhat_1)(s)
\; ds\Big\|_{L^1}
\\
\notag
&\qquad +
\Big\| \intl_0^t 
  \TT_{3,\epsi}(s) \T(\uhat_1,D_\mu \uhat_1,\uhat_1)(s)
\; ds\Big\|_{L^1}
\\
\label{Proof.Prop02.Eq01}
&\qquad +
\Big\| \intl_0^t 
  \TT_{3,\epsi}(s) \T(\uhat_1,\uhat_1,D_\mu \uhat_1)(s)
\; ds\Big\|_{L^1}.
\end{align}
Each of the three terms in \eqref{Proof.Prop02.Eq01} can be treated by adapting the estimates from Case~2 of the proof of Proposition~\ref{Proposition01}. 
Instead of Lemma~\ref{Lemma.bound.uhat1dot} one has to use that
\begin{align*}
\sup_{s\in[0,\tstar/\epsi]} \| \pt D_\mu \P_\epsi \uhat_1(s) \|_{L^1}
\end{align*}
is uniformly bounded, which can be shown in a similar way. This leads to \eqref{Gronwall.bound} for $j=3$ with constants $C_1$ and $C_2$ which depend on $ C_{u,2}$.
For $j=1$ the bound \eqref{Gronwall.bound} can be shown by adapting the procedure from Section~\ref{Subsubsec.proof.prop01.term.1}.
\qed

\section{Error bound for the approximation}\label{Sec.error.bound.approximation}

If $ u_1, u_3 $ is the solution of \eqref{PDE.mfe} with initial data \eqref{PDE.mfe.initial.data}
then $\uutilde$ defined in \eqref{Ansatz}
provides an approximation to the exact solution $ \uu $ of the original problem \eqref{PDE.uu}. 
Our goal is now to prove an error bound for this approximation.
This error bound requires an additional assumption.

\begin{Assumption}[Non-resonance condition]\label{Ass:Delta.lambda.j5}
The matrices $\L_3(0) = \L(3\w, 3\kappa)$ and $\L_5(0) $ $ = \L(5\w, 5\kappa)$ have no common eigenvalues, i.e.
$ \lambda_{3\ell}(0) \not= \lambda_{5m}(0) $ for all $\ell,m=1, \ldots, n$.
\end{Assumption}

Local well-posedness of \eqref{PDE.uu} in the Wiener algebra on long time intervals $ [0,\tend/\epsi] $ for some $\tend>0$ could be shown by adapting the proof of Lemma~\ref{Lemma.wellposedness}. We may thus assume that a unique mild solution of 
\eqref{PDE.uu} exists on $[0,\tstar/\epsi]$ (after decreasing the value of $\tstar$ if necessary).

\mynote{Der Operator $ A(\partial) + \tfrac{1}{\epsi}E $ hat die gleichen Eigenschaften wie der Operator $\A$ aus dem Beweis von Lemma~\ref{Lemma.wellposedness}. Insbesondere erhalten die Gruppenoperatoren auf $W$ die Norm. Die Algebra-Eigenschaft braucht man nat\"urlich auch. Dass $\uu$ im Raum oszilliert, spielt keine Rolle, da man ja im Wohlgestelltheitsbeweis nie 
$ \big(A(\partial) + \tfrac{1}{\epsi}E\big)\uu $ betrachtet, sondern nur $ \exp\big(t (A(\partial) + \tfrac{1}{\epsi}E)\big)\uu $.
Damit kann man das Fixpunktargument genauso durchziehen wie im Lemma.
}

\begin{Theorem}\label{Theorem.error.bound}
Let $ p \in W^2 $ and let $\uu$ be the solution of \eqref{PDE.uu}.
Let $(u_1,u_3)$ be the classical solution of \eqref{PDE.mfe}
established in part~\ref{Lemma.wellposedness.3}
of Lemma~\ref{Lemma.wellposedness}, and let $\uutilde$ be the approximation defined in \eqref{Ansatz}.
Under Assumptions~\ref{Ass:polarization}, \ref{Ass:L.properties}, and \ref{Ass:Delta.lambda.j5} there is a constant such that 
\begin{align}
\label{Theorem.error.bound.assertion01}
 \sup_{t\in[0,\tstar/\epsi]} \| \uu(t) - \uutilde(t) \|_W
 &\leq C\epsi^2,
 \\
\label{Theorem.error.bound.assertion02}
 \sup_{t\in[0,\tstar/\epsi]} \| \uu(t) - \uutilde(t) \|_{L^\infty}
 &\leq C\epsi^2.
\end{align}
\end{Theorem}
\bigskip

Note that Proposition~\ref{Proposition01} and \ref{Proposition02} apply
under the assumptions of Theorem~\ref{Theorem.error.bound}, which yields the bounds
\eqref{Bounds.uhat}, \eqref{Bounds.z}, \eqref{Bounds.Dmu.uhat}, \eqref{Bounds.Dmu.z} in addition to 
\eqref{Lemma.wellposedness.bound.2}.
Applying \eqref{Bounds.uhat} and \eqref{Bounds.Dmu.uhat} to 
\begin{align*}
\| u_3(t) \|_{W^1} &= \| u_3(t) \|_W + \sum_{\mu=1}^d \| \partial_\mu u_3(t) \|_W
= \| \uhat_3(t) \|_{L^1} + \sum_{\mu=1}^d \| D_\mu \uhat_3(t) \|_{L^1} 
\end{align*}
yields in particular that
\begin{align}
\label{Bound.u3.W1}
\sup_{t\in[0,\tstar/\epsi]} \| u_3(t) \|_{W^1} \leq C\epsi.
\end{align}
Moreover, \eqref{zdot.in.terms.of.uhat}, \eqref{Bounds.z} and Lemma~\ref{Lemma.useful.bounds} yield
\begin{align}
\sup_{t\in[0,\tstar/\epsi]} \| \pt z_{j}(t) \|_{L^1}
&\leq
\epsi \sup_{t\in[0,\tstar/\epsi]} \sum_{\#J=j} \| \FF(t,\uhat,J)(t) \|_{L^1}
\label{Bound.zdot}
\leq
C \sum_{\#J=j} \epsi^{1+(|J|_1-3)/2} \leq C\epsi.
\end{align}

\noindent
\textbf{Proof.}
Since the proof of Theorem~\ref{Theorem.error.bound} is rather lengthy, we subdivide it into several steps.

\paragraph{Step 1.} Our first goal is to derive an evolution equation for the error 
$ \delta=\uu-\uutilde $ and its Fourier transform.
The approximation $\uutilde$ solves \eqref{PDE.uu} up to the residual
\begin{align}
\label{Def.R}
R(t,x) =
\epsi T(\uutilde,\uutilde,\uutilde)(t,x) - 
\Big(\pt \uutilde(t,x) +A(\partial)\uutilde(t,x) +\frac{1}{\epsi} E\uutilde(t,x)\Big).
\end{align}
In order to derive a more useful expression for $R$, we note that \eqref{PDE.mfe} yields
\begin{align}
\nonumber
& \pt \uutilde(t,x) +A(\partial)\uutilde(t,x) +\frac{1}{\epsi} E\uutilde(t,x)
\\
\nonumber
&=
 \sum_{j\in\J} 
 e^{\ii j (\kappa \cdot x - \w t)/\epsi}
\Big(
 \pt u_j(t,x) + \tfrac{\ii}{\epsi}\L(j\w,j\kappa)u_j(t,x) + A(\partial)u_j(t,x) 
\Big) 
 \\
 &=
 \epsi
 \sum_{j\in\J} 
 \sum_{\#J=j} 
 e^{\ii j (\kappa \cdot x - \w t)/\epsi}
 T(u_{j_1},u_{j_2},u_{j_3})(t,x)
 \label{R.sum.1}
 \end{align}
in contrast to
\begin{align}
\nonumber
\epsi T(\uutilde,\uutilde,\uutilde)(t,x)
&= 
\epsi 
\sum_{J\in\J^3} 
 e^{\ii \#J (\kappa \cdot x - \w t)/\epsi}
T(u_{j_1},u_{j_2},u_{j_3})(t,x)
\\
&= 
\epsi
\sum_{\substack{j \text{ odd} \\ |j|\leq 9}}
\sum_{\#J=j} 
 e^{\ii j (\kappa \cdot x - \w t)/\epsi}
T(u_{j_1},u_{j_2},u_{j_3})(t,x).
\label{R.sum.2}
\end{align}
The difference is that \eqref{R.sum.2} includes summands with $|j| \in\{5,7,9\}$, 
whereas $ j\in\J $ implies that 
$ |j|\in\{1,3\}$ in \eqref{R.sum.1}.
Inserting \eqref{R.sum.1} and \eqref{R.sum.2} into \eqref{Def.R} yields
\begin{align*}
R(t,x) &= \epsi \sum_{|j| \in\{5,7,9\}}
\sum_{\#J=j} 
 e^{\ii j (\kappa \cdot x - \w t)/\epsi}
 T(u_{j_1},u_{j_2},u_{j_3})(t,x),
\end{align*}
and by definition of $R$ the error $ \delta=\uu-\uutilde $ solves the equation
\begin{align}
 \label{PDE.delta}
\pt \delta =  
- A(\partial)\delta
- \frac{1}{\epsi} E\delta
+ \epsi\left[T(\uu,\uu,\uu)-T(\uutilde,\uutilde,\uutilde)\right]
+ R.
\end{align}
In order to derive a bound in $\| \cdot \|_W$, we apply the Fourier transform to
\eqref{PDE.delta} to obtain
\begin{align*}
\pt \widehat{\delta}(t,k) =  
-\big(\ii A(k)+\tfrac{1}{\epsi} E\big)
\widehat{\delta}(t,k) 
+ \epsi \G\big(\F\uu,\F\uutilde \big)(t,k)
+ \widehat{R}(t,k) 
\end{align*}
with 
\begin{align}
\nonumber
\G\big(\F\uu,\F\uutilde\big)&=
\T(\F\uu,\F\uu,\F\uu)-\T(\F\uutilde,\F\uutilde,\F\uutilde),
\\
\nonumber
\widehat{R}(t,k) 
&= 
\epsi
\sum_{|j| \in\{5,7,9\}} \sum_{\#J=j} 
\F\left(
 T(u_{j_1},u_{j_2},u_{j_3})
 e^{\ii j \kappa \cdot x /\epsi} \right)(t,k) e^{-\ii j \w t/\epsi} 
 \\
 \label{Def.Rhat}
 &= 
 \epsi 
 \sum_{|j| \in\{5,7,9\}}  \sum_{\#J=j} 
 \T(\uhat_{j_1},\uhat_{j_2},\uhat_{j_3})(t,k-\tfrac{j\kappa}{\epsi})
e^{-\ii j \w t/\epsi}.  
\end{align}

\paragraph{Step 2.} 
With the variation-of-constants formula $\widehat{\delta}(t,k) $
can be expressed as
\begin{align}
\nonumber
\widehat{\delta}(t,k) 
&=  
\epsi \intl_0^t
\exp\big((s-t)\big(\ii A(k)+\tfrac{1}{\epsi} E\big)\big)
\G\big(\F\uu(s),\F\uutilde(s) \big)(k) \; ds
\\[-3mm]
\label{Theorem.error.bound.01}
\\[-3mm]
\nonumber
& \quad +
\intl_0^t
\exp\big((s-t)\big(\ii A(k)+\tfrac{1}{\epsi} E\big)\big)
\widehat{R}(s,k) \; ds.
\end{align}
We aim for proving 
$ \sup_{t\in[0,\tstar/\epsi]} \| \widehat{\delta}(t) \|_{L^1} \leq C \epsi^2 $
via Gronwall's lemma.
Since the matrix $\ii A(k)+\tfrac{1}{\epsi} E $ is skew-hermitian for every $k$,
the first term on the right-hand side of \eqref{Theorem.error.bound.01} can be bounded in $L^1$ by
\begin{align}
\nonumber
& \epsi \intl_0^t
 \intl_{\IR^d}
 \left| \exp\big((s-t)\big(\ii A(k)+\tfrac{1}{\epsi} E\big)\big) \right|_2
\big| \G\big(\F\uu(s),\F\uutilde(s) \big)(k) \big|_2 \; dk \; ds
\\
\nonumber
&\leq 
 \epsi \intl_0^t
 \intl_{\IR^d}
| \G\big(\F\uu(s),\F\uutilde(s) \big)(k) |_2 \; dk \; ds
\\
&\leq 
3 C_\T C_\uu^2 \; 
\epsi \intl_0^t
\| \widehat{\delta}(s) \|_{L^1} \; ds
\label{Theorem.error.bound.02}
\end{align}
due to \eqref{Bound.T.2} where $C_\uu$ is a constant such that
\begin{align*}
 \sup_{t\in[0,\tstar/\epsi]} \| \uu(t) \|_W \leq C_\uu, \qquad
 \sup_{t\in[0,\tstar/\epsi]} \| \uutilde(t) \|_W \leq C_\uu
\end{align*}
uniformly in $\epsi$.
The laborious part of the proof is to show that
\begin{align}
 \label{Theorem.error.bound.crucial}
 \sup_{t\in[0,\tstar/\epsi]}
 \Big\|
 \intl_0^t
\exp\big((s-t)\big(\ii A(\cdot)+\tfrac{1}{\epsi} E\big)\big)
\widehat{R}(s) \; ds
\Big\|_{L^1} \leq C \epsi^2
\end{align}
with a constant $C$ which does not depend on $\epsi$.
If \eqref{Theorem.error.bound.crucial} holds, then together with
\eqref{Theorem.error.bound.01} and \eqref{Theorem.error.bound.02}
it follows that
\begin{align*}
\| \widehat{\delta}(t) \|_{L^1} 
&\leq 
C C_\uu^2 \; \epsi \intl_0^t
\| \widehat{\delta}(s) \|_{L^1} \; ds
+ C \epsi^2,
\end{align*}
and Gronwall's lemma yields 
\begin{align*}
\sup_{t\in[0,\tstar/\epsi]} \| \uu(t)-\uutilde(t) \|_W 
=
\sup_{t\in[0,\tstar/\epsi]} \| \widehat{\delta}(t) \|_{L^1} 
&\leq 
C \epsi^2 e^{\gamma \tstar}
\end{align*}
with $\gamma=C C_\uu^2$, which proves \eqref{Theorem.error.bound.assertion01}.
The second bound \eqref{Theorem.error.bound.assertion02}
is an immediate consequence of the embedding $W(\IR^d) \hookrightarrow L^\infty(\IR^d)$.
The purpose of the following steps is to prove \eqref{Theorem.error.bound.crucial}.

\paragraph{Step 3.}
In this step the term on the left-hand side of \eqref{Theorem.error.bound.crucial} is reformulated. With
the change of variables $k' = k-\tfrac{j\kappa}{\epsi}$, 
$\epsi k = j \kappa + \epsi k' $, definition \eqref{Def.L.j} and the representation \eqref{Def.Rhat}, we obtain
\begin{align*}
& \intl_0^t
\exp\big((s-t)\big(\ii A(k)+\tfrac{1}{\epsi} E\big)\big)
\widehat{R}(s,k) \; ds
\\
&=
 \epsi 
 \sum_{|j| \in\{5,7,9\}}  \sum_{\#J=j} 
\intl_0^t
\exp\left(\frac{\ii}{\epsi}(s-t)\L(j\w,\epsi k)\right) 
e^{-\ii j \w t/\epsi}
\T(\uhat_{j_1},\uhat_{j_2},\uhat_{j_3})(s,k-\tfrac{j\kappa}{\epsi})
 \; ds
\\
&=
 \epsi e^{-\ii j \w t/\epsi}
 \sum_{|j| \in\{5,7,9\}}  \sum_{\#J=j} 
\intl_0^t
\exp\left(\frac{\ii}{\epsi}(s-t)\L_j(\epsi k')\right) 
 \T(\uhat_{j_1},\uhat_{j_2},\uhat_{j_3})(s,k')
 \; ds.
\end{align*}

\mynote{
\begin{align*}
& \intl_0^t
\exp\big((s-t)\big(\ii A(k)+\tfrac{1}{\epsi} E\big)\big)
\widehat{R}(s) \; ds
\\
&=
 \epsi 
\intl_0^t
\exp\big((s-t)\big(\ii A(k)+\tfrac{1}{\epsi} E\big)\big) 
\\
&\hs{20} \times
 \sum_{|j| \in\{5,7,9\}}  \sum_{\#J=j} 
 \T(\uhat_{j_1},\uhat_{j_2},\uhat_{j_3})(s,k-\tfrac{j\kappa}{\epsi})
e^{-\ii j \w s/\epsi} \; ds
\\
&=
 \epsi 
 \sum_{|j| \in\{5,7,9\}}  \sum_{\#J=j} 
\intl_0^t
\exp\left(\frac{\ii}{\epsi}(s-t)\big(-j \w + A(\epsi k)-\ii E\big)\right) 
e^{-\ii j \w t/\epsi}
\\
&\hs{30} \times
 \T(\uhat_{j_1},\uhat_{j_2},\uhat_{j_3})(s,k-\tfrac{j\kappa}{\epsi})
 \; ds
\\
&=
 \epsi 
 \sum_{|j| \in\{5,7,9\}}  \sum_{\#J=j} 
\intl_0^t
\exp\left(\frac{\ii}{\epsi}(s-t)\L(j\w,\epsi k)\right) 
e^{-\ii j \w t/\epsi}
\\
&\hs{30} \times
 \T(\uhat_{j_1},\uhat_{j_2},\uhat_{j_3})(s,k-\tfrac{j\kappa}{\epsi})
 \; ds.
\end{align*}
}

From now on, we omit the dash and write $k$ instead of $k'$. The difference does not matter because later we integrate over $k$, anyway.
By \eqref{Def.TT} and \eqref{Def.F.uhat} we have
\begin{align*}
\exp\left(\frac{\ii}{\epsi}(s-t)\L_j(\epsi k)\right) \T(\uhat_{j_1},\uhat_{j_2},\uhat_{j_3})(s,k)
&=
\exp\left(-\frac{\ii t}{\epsi}\L_j(\epsi k)\right) \Psi_{j}(\epsi k) \FF(s,\uhat,J)(k),
\end{align*}
which yields the bound
\begin{align}
\label{Theorem.error.bound.03}
 \Big\|
 \intl_0^t
\exp\big((s-t)\big(\ii A(\cdot)+\tfrac{1}{\epsi} E\big)\big)
\widehat{R}(s) \; ds
\Big\|_{L^1} 
&\leq
 \epsi 
 \sum_{|j| \in\{5,7,9\}}  \sum_{\#J=j} 
\Big\| \intl_0^t \FF(s,\uhat,J)
\; ds \Big\|_{L^1}
\end{align}
for the crucial term. 

\paragraph{Step 4.}

Our next goal is to prove that
\begin{align*}
 \sum_{|j| \in\{5,7,9\}}  \sum_{\#J=j} 
\Big\| \intl_0^t \FF(s,\uhat,J)
\; ds \Big\|_{L^1} \leq C \epsi.
\end{align*}
Combining this estimate with \eqref{Theorem.error.bound.03} yields the desired bound \eqref{Theorem.error.bound.crucial}.
Note that there is an extra factor $\epsi$ on the right-hand side of \eqref{Theorem.error.bound.03}.

Let $j$ be an odd number with $|j| \in\{5,7,9\}$, and let $J\in\J^3$ with $\#J=j$. 
We distinguish the two cases 
$|J|_1=|j|=5$ and $|J|_1 \in \{7,9\}$. The latter is the easier one, because in this case 
Lemma~\ref{Lemma.useful.bounds}, \eqref{Lemma.wellposedness.bound.2}, \eqref{Norm.uhat.z}, and \eqref{Bounds.uhat} yield
\begin{align*}
\Big\| \intl_0^t \FF(s,\uhat,J) \; ds \Big\|_{L^1} 
&\leq
\intl_0^t \| \FF(s,\uhat,J) \|_{L^1} \; ds
\\
& \leq
C_\T \frac{\tstar}{\epsi} \sup_{s\in [0,\tstar/\epsi]} 
\| \uhat_{j_1}(s) \|_{L^1} \| \uhat_{j_2}(s) \|_{L^1} \| \uhat_{j_3}(s) \|_{L^1} 
\\
& \leq
C\tstar \epsi^{(|J|_1-5)/2}
\\
& \leq
C \tstar \epsi
\end{align*}
for $|J|_1 \in \{7,9\}$.
Now let $|J|_1=|j|=5$. This situation appears only if $ J$ is either $(3,1,1)$, $(1,3,1)$, or $(1,1,3)$.
Since all three cases can be treated in the same way, we may henceforth assume that $J=(3,1,1)$, which means that 
\begin{align*}
\FF(t,\uhat,J) = \TT_{5,\epsi}(t) \T\big(\uhat_3,\uhat_1,\uhat_1\big)(t)
\end{align*}
by \eqref{Def.F.uhat}.
As in Case 2 in Section~\ref{Subsubsec.proof.prop01.term.3} we decompose
\begin{align*}
\T\big(\uhat_3,\uhat_1,\uhat_1\big)
&=
\T\big(\uhat_3,\P_\epsi\uhat_1,\P_\epsi\uhat_1\big)
+ \T\big(\uhat_3,\P_\epsi\uhat_1,\P_\epsi^\perp\uhat_1\big)
\\
&\quad
+ \T\big(\uhat_3,\P_\epsi^\perp\uhat_1,\P_\epsi\uhat_1\big)
+ \T\big(\uhat_3,\P_\epsi^\perp\uhat_1,\P_\epsi^\perp\uhat_1\big),
\end{align*}
and as before all terms involving $\P_\epsi^\perp\uhat_1$ can be treated in a straightforward way because
of \eqref{Bounds.uhat}.
Hence, the main difficulty is to prove that
\begin{align}
\label{Main.difficulty}
\Big\| \intl_0^t \TT_{5,\epsi}(s) \T\big(\uhat_3,\P_\epsi\uhat_1,\P_\epsi\uhat_1 \big)(s) \; ds \Big\|_{L^1}
\leq C \epsi.
\end{align}

\paragraph{Step 5.}

To prove \eqref{Main.difficulty}, we use that \eqref{Def.z} and \eqref{Def.TT} yield the representation
\begin{align}
\uhat_3(t,k) &= \TT_{3,\epsi}^*(t,k) z_3(t,k)
\notag
\\
\notag
&= 
\Psi_3(\epsi k) \exp\big(-\tfrac{\ii t}{\epsi}\Lambda_3(\epsi k)\big) z_3(t,k) 
\\
&=
\label{representation.uhat.3}
\sum_{m=1}^n 
\exp\big(-\tfrac{\ii t}{\epsi}\lambda_{3m}(\epsi k)\big)z_{3m}(t,k)\psi_{3m}(\epsi k).
\end{align}
Recall that $\psi_{jm} \in \IC^n $ is the $m$-th column of the unitary matrix $ \Psi_j $, and that
$\lambda_{jm} $ is the $m$-th entry on the diagonal of $ \Lambda_j \in\IR^{n \times n} $; 
cf.~Section~\ref{Subsec.Transformation}.
A similar representation of $ \P_\epsi\uhat_1 $ is derived by using \eqref{Relation.between.projections}
and \eqref{Def.P1} in addition to \eqref{Def.z} and \eqref{Def.TT}, namely
\begin{align}
\P_\epsi\uhat_1(t,k) &= \TT_{1,\epsi}^*(t,k) Pz_1(t,k)
\notag
\\
\notag
&= 
\Psi_1(\epsi k) \exp\big(-\tfrac{\ii t}{\epsi}\Lambda_1(\epsi k)\big) Pz_1(t,k) 
\\
&=
\label{representation.Puhat.1}
\exp\big(-\tfrac{\ii t}{\epsi}\lambda_{11}(\epsi k)\big)z_{11}(t,k)\psi_{11}(\epsi k).
\end{align}
Combining the representations \eqref{representation.uhat.3} and \eqref{representation.Puhat.1} 
with \eqref{Def.T.Fourier} and using $ \Psi_{5}^* = \sum_{\ell=1}^n e_\ell \psi_{5\ell}^* $ yields
\begin{align}
\label{TTT.F}
\TT_{5,\epsi}(s) \T\big(\uhat_3,\P_\epsi\uhat_1,\P_\epsi\uhat_1 \big)(s)
=
F(s,z)
\end{align}
with $F$ defined by
\begin{align*}
F(t,z) &= \Big(F_{\ell}(t,z)\Big)_{\ell=1}^n,
\\
F_{\ell}(t,z)(k)
&=
\sum_{m=1}^n \intl_{\#K = k} 
\exp\Big(\tfrac{\ii t}{\epsi} \Delta\lambda_{\ell m}(\epsi,k,K) \Big)
Z_m(t,K) c_{\ell m}(\epsi,k,K) \; dK
\end{align*}
and with the notation 
\begin{align*}
 K&=\big(k^{(1)},k^{(2)},k^{(3)}\big)\in \IR^d \times \IR^d \times \IR^d, 
 \\
\Delta\lambda_{\ell m}(\epsi,k,K) &= \lambda_{5\ell}(\epsi k)-\lambda_{3m}\big(\epsi k^{(1)}\big) 
- \lambda_{11}\big(\epsi k^{(2)}\big) - \lambda_{11}\big(\epsi k^{(3)}\big),
\\
Z_m(t,K) &= z_{3m}\big(t,k^{(1)}\big) z_{11}\big(t,k^{(2)}\big) z_{11}\big(t,k^{(3)}\big),
\\
c_{\ell m}(\epsi,k,K) &=
\frac{1}{(2\pi)^d} \psi_{5\ell}^*(\epsi k)
T\Big(\psi_{3m}\big(\epsi k^{(1)}\big),\psi_{11}\big(\epsi k^{(2)}\big),\psi_{11}\big(\epsi k^{(3)}\big)\Big).
\end{align*}
Since by definition $ \psi_{j\ell} $ is the $\ell$--th column of the unitary matrix $ \Psi_j $, it follows that
\begin{align*}
| c_{\ell m}(\epsi,k,K) | \leq C_\T \qquad \text{for all } \epsi, k, K.
\end{align*}
With \eqref{TTT.F} the left-hand side of \eqref{Main.difficulty} can be bounded by
\begin{align}
\notag
& \Big\| \intl_0^t \TT_{5,\epsi}(s) \T\big(\uhat_3,\P_\epsi\uhat_1,\P_\epsi\uhat_1 \big)(s) \; ds \Big\|_{L^1}
\\
\notag
&=
 \intl_{\IR^d} \Big| \intl_0^t F(s,z)(k) \; ds \Big|_2 \; dk
\\
\notag
&\leq
 \sum_{\ell=1}^n \intl_{\IR^d} \Big| \intl_0^t F_\ell(s,z)(k) \; ds \Big| \; dk
\\
\notag
&\leq
 \sum_{\ell,m=1}^n \intl_{\IR^d} \Big| 
\intl_{\#K = k} \intl_0^t 
\exp\Big(\tfrac{\ii s}{\epsi} \Delta\lambda_{\ell m}(\epsi,k,K) \Big)
Z_m(s,K) \; ds  \; c_{\ell m}(\epsi,k,K) \; dK \Big| \; dk
\\
\notag
&\leq
 \sum_{\ell,m=1}^n \intl_{\IR^d}  
\intl_{\#K = k} \Big|
\intl_0^t \exp\Big(\tfrac{\ii s}{\epsi} \Delta\lambda_{\ell m}(\epsi,k,K) \Big)
Z_m(s,K) \; ds  \Big|
\; | c_{\ell m}(\epsi,k,K) |  \; dK \; dk
\\
\label{Theorem.error.bound.04}
&\leq
C_\T
 \sum_{\ell,m=1}^n \intl_{\IR^d}  
\intl_{\#K = k} \Big|
\intl_0^t \exp\Big(\tfrac{\ii s}{\epsi} \Delta\lambda_{\ell m}(\epsi,k,K) \Big)
Z_m(s,K) \; ds  \Big|
\;  \; dK \; dk.
\end{align}

\paragraph{Step 6.}

In this step, we prove that 
\begin{align}
\label{Goal.step.6}
&
\Big| \intl_0^t \exp\Big(\tfrac{\ii s}{\epsi} \Delta\lambda_{\ell m}(\epsi,k,K) \Big)
Z_m(s,K) \; ds  \Big|
\\
\nonumber
&\leq
C \epsi \Big(|Z_m(t,K)| + \sum_{i=1}^3  |k^{(i)} |_1 \intl_0^t |Z_m(s,K)| \; ds
+ \intl_0^t | \pt Z_m(s,K) | \; ds \Big)
\end{align}
for all $\ell, m, K$ and $k=\#K$. Since $\ell, m, k, K$ are considered \emph{fixed} in this step, we can 
simplify notation by setting
\begin{align*}
\Delta\lambda(\epsi) &= \Delta\lambda_{\ell m}(\epsi,k,K),
\\
Z(s) &= Z_m(s,K), 
\\
Y(s) &= \exp\Big(\tfrac{\ii s}{\epsi} \big[\Delta\lambda(\epsi) - \Delta\lambda(0) \big] \Big) Z(s).
\end{align*}
Since $\lambda_{11}(0)=0$, we have
\begin{align*}
\Delta\lambda(0)  =  \lambda_{5\ell}(0)-\lambda_{3m}(0) \not= 0
\end{align*}
by Assumption~\ref{Ass:Delta.lambda.j5}. Hence, we can integrate by parts in order to generate one additional power of $\epsi$.
Now we obtain for the left-hand side of \eqref{Goal.step.6}
\begin{align*}
& \Big|
\intl_0^t \exp\Big(\tfrac{\ii s}{\epsi} \Delta\lambda(\epsi) \Big) Z(s) \; ds  
\Big|
=
\Big|
\intl_0^t \exp\Big(\tfrac{\ii s}{\epsi} \Delta\lambda(0) \Big) Y(s) \; ds  
\Big|
\\
&\leq
\frac{\epsi}{|\Delta\lambda(0)|} \Big| \Big[  
 \exp\Big(\tfrac{\ii s}{\epsi} \Delta\lambda(0) \Big) Y(s) \Big]_{s=0}^t
 \Big|
+
 \frac{\epsi}{|\Delta\lambda(0)|} \Big| \intl_0^t 
 \exp\Big(\tfrac{\ii s}{\epsi} \Delta\lambda(0) \Big) 
 \pt Y(s) \; ds \Big|
 \\
 &\leq
C \epsi \Big(|Z(t)|+|Z(0)|\Big)
+
C \epsi \intl_0^t | \pt Y(s) | \; ds. 
\end{align*}
The initial conditions \eqref{PDE.mfe.Fourier.initial.data} imply that $ z_{3m}(0,k)=0 $ for all $m$ and hence that
$ Z(0)=Z_m(0,K)=0$.
For the last term, we obtain
\begin{align}
\label{Irgendein.label}
\intl_0^t | \pt Y(s) | \; ds 
\leq \frac{1}{\epsi} \intl_0^t 
\Big| \big[\Delta\lambda(\epsi) - \Delta\lambda(0) \big] \Big| \cdot | Z(s) | \; ds 
+
\intl_0^t 
| \pt Z(s) | \; ds.
\end{align}
For the difference
\begin{align*}
\Delta\lambda(\epsi) - \Delta\lambda(0)
&=
\Delta\lambda_{\ell m}(\epsi) - \Delta\lambda_{\ell m}(0)
\\
&=
\big(\lambda_{5\ell}(\epsi k) - \lambda_{5\ell}(0) \big)
- \big(\lambda_{3m}\big(\epsi k^{(1)}\big) - \lambda_{3m}(0)\big)
\\
& \quad
- \big(\lambda_{11}\big(\epsi k^{(2)}\big) - \lambda_{11}(0)\big)
- \big(\lambda_{11}\big(\epsi k^{(3)}\big) - \lambda_{11}(0)\big)
\end{align*}
the Lipschitz continuity of the eigenvalues (see Assumption~\ref{Ass:L.properties}\ref{Ass:L.properties.item.ii}) and the identity
$ k= \#K = k^{(1)}+k^{(2)}+k^{(3)} $ yield the inequality
\begin{align*}
| \Delta\lambda(\epsi) - \Delta\lambda(0) | 
\leq
C |\epsi k|_1 + C \sum_{i=1}^3 |\epsi k^{(i)}|_1
\leq
2C\epsi \sum_{i=1}^3 |k^{(i)}|_1,
\end{align*}
and the $\epsi$ on the right-hand side cancels with the factor $1/\epsi$ in \eqref{Irgendein.label}.
All in all, this yields the bound
\begin{align*}
\Big| \intl_0^t \exp\Big(\tfrac{\ii s}{\epsi} \Delta\lambda(\epsi) \Big)
Z(s) \; ds  \Big|
&\leq 
C \epsi \Big(|Z(t)| + \sum_{i=1}^3  |k^{(i)} |_1 \intl_0^t | Z(s)| \; ds
+ \intl_0^t | \pt Z(s) | \; ds \Big),
\end{align*}
which proves \eqref{Goal.step.6}.

\paragraph{Step 7.}

Now we make our way back to \eqref{Main.difficulty}. 
Substituting \eqref{Goal.step.6} into \eqref{Theorem.error.bound.04} yields
\begin{align*}
& \Big\| \intl_0^t \TT_{5,\epsi}(s) \T\big(\uhat_3,\P_\epsi\uhat_1,\P_\epsi\uhat_1 \big)(s) \; ds \Big\|_{L^1}
\\
&\leq
C_\T
\sum_{\ell,m=1}^n \intl_{\IR^d}  
\intl_{\#K = k} \Big|
\intl_0^t \exp\Big(\tfrac{\ii s}{\epsi} \Delta\lambda_{\ell m}(\epsi,k,K) \Big)
Z_m(s,K) \; ds  \Big|
\;  \; dK \; dk
\\
&\leq
C_\T n \epsi
\Big(X_1(t)+X_2(t)+X_3(t)\Big)
\end{align*}
with the short-hand notation
\begin{align*}
X_1(t) &=
\sum_{m=1}^n \intl_{\IR^d}  
\intl_{\#K = k} |Z_m(t,K)| \;  \; dK \; dk,
\\
X_2(t) &=
\sum_{i=1}^3 \sum_{m=1}^n \intl_{\IR^d}  
\intl_{\#K = k} |k^{(i)} |_1 \intl_0^t |Z_m(s,K)| \; ds
\;  \; dK \; dk,
\\
X_3(t) &=
\sum_{m=1}^n \intl_{\IR^d}  
\intl_{\#K = k} \intl_0^t | \pt Z_m(s,K) | \; ds 
\;  \; dK \; dk.
\end{align*}
To conclude the proof, it has to be shown that $ X_i(t) \leq C $ for $i=1,2,3$ and for all $t\in [0,\tstar/\epsi]$ with a constant $c$ that does not depend on $\epsi$.
From the definition of $Z_m$ it follows that 
\begin{align*}
\sum_{m=1}^n  |Z_m(t,K)| 
&=
\Big(\sum_{m=1}^n \big| z_{3m}\big(t,k^{(1)}\big) \big|\Big) \cdot \big| z_{11}\big(t,k^{(2)}\big) \big| \cdot \big| z_{11}\big(t,k^{(3)}\big) \big|
\\
&=
\big| z_3\big(t,k^{(1)}\big) \big|_1 \cdot \big| z_{11}\big(t,k^{(2)}\big) \big| \cdot \big| z_{11}\big(t,k^{(3)}\big) \big|
\\
&\leq
\sqrt{n} 
\big| z_3\big(t,k^{(1)}\big) \big|_2 \cdot \big| Pz_1\big(t,k^{(2)}\big) \big|_2 \cdot \big| Pz_1\big(t,k^{(3)}\big) \big|_2
\\
&\leq
\sqrt{n} 
\big| \uhat_3\big(t,k^{(1)}\big) \big|_2 \cdot \big| \uhat_1\big(t,k^{(2)}\big) \big|_2 \cdot \big| \uhat_1\big(t,k^{(3)}\big) \big|_2
\end{align*}
due to \eqref{Norm.uhat.z}. This yields
\begin{align*}
X_1(t) &=
\sum_{m=1}^n \intl_{\IR^d} \intl_{\#K = k} |Z_m(t,K)| \;  \; dK \; dk
\\
&\leq
\sqrt{n} \intl_{\IR^d} \intl_{\#K = k} 
\big| \uhat_3\big(t,k^{(1)}\big) \big|_2 \cdot \big| \uhat_1\big(t,k^{(2)}\big) \big|_2 \cdot \big| \uhat_1\big(t,k^{(3)}\big) \big|_2
\;  \; dK \; dk
\\
&= 
\sqrt{n} 
\| u_3(t) \big\|_W \cdot \| u_1(t) \|_W^2 
\leq 
C.
\end{align*}
(In fact, we even have that $X_1(t) \leq C\epsi$ according to \eqref{Bound.uhat.3}.)
To bound $X_2(t)$ we use that
\begin{align*}
\intl_{\IR^d}
|k^{(i)}|_1 \big| \uhat_{j_i}\big(t,k^{(i)}\big) \big|_2
\; dk^{(i)} \leq \| u_{j_i}(t) \big\|_{W^1}
\end{align*}
for $ i\in\{1,2,3\}$ and $J=(j_1,j_2,j_3)=(3,1,1)$ to obtain in a similar way as before that
\begin{align*}
X_2(t) &=
\sum_{i=1}^3 \sum_{m=1}^n \intl_{\IR^d} \intl_{\#K = k}
|k^{(i)} |_1 \intl_0^t |Z_m(s,K)| \; ds
\;  \; dK \; dk
\\
&\leq
C \intl_0^t 
\| u_3(s) \big\|_{W^1} \cdot \| u_1(s) \|_{W^1}^2
\; ds
\\
&\leq
C \tfrac{\tstar}{\epsi}
\sup_{s \in [0,\tstar/\epsi]} \big(\| u_3(s) \|_{W^1} \cdot \| u_1(s) \|_{W^1}^2\big)
\leq C
\end{align*}
due to \eqref{Bound.u3.W1} and \eqref{Lemma.wellposedness.bound.2}.
Finally, it follows from \eqref{Bounds.z} and \eqref{Bound.zdot} that
\begin{align*}
X_3(t) &= \sum_{m=1}^n \intl_{\IR^d}  \intl_{\#K = k}
\intl_0^t | \pt Z_m(s,K) | \; ds 
\;  \; dK \; dk
\\
&\leq
C \intl_0^t 
\Big(
\| \pt z_3(s) \big\|_{L^1} \cdot \| z_1(s) \|_{L^1}^2
\\
& \hs{15}
+ 2 \| z_3(s) \big\|_{L^1} \cdot \| z_1(s) \|_{L^1}  \cdot \| \pt z_1(s) \|_{L^1} 
\Big)
\; ds
\\
&\leq C\tfrac{\tstar}{\epsi}(\epsi+2 \epsi^2) \leq C.
\end{align*}
This completes the proof of Theorem~\ref{Theorem.error.bound}.
\qed

\section{Discussion}\label{Sec.higher.order.approximation}

In this work the solution $\uu$ of the semilinear PDE system \eqref{PDE.uu} is approximated by
\begin{align}
\label{Ansatz.discussion}
\uutilde(t,x) = \sum_{j\in \J}
e^{\ii j(\kappa \cdot x - \w t)/\epsi} u_j(t,x),
\qquad u_{-j} = \overline{u_j}, \qquad \J = \{\pm 1, \pm 3\}
\end{align}
with $ u_j$ solving \eqref{PDE.mfe}--\eqref{PDE.mfe.initial.data}.
Theorem~\ref{Theorem.error.bound} states that the error of the approximation on the long time interval
$ [0,\tstar/\epsi] $ is $\Ord{\epsi^2}$.

A natural question to ask is whether or not the approximation improves if $\J$ is replaced by
$\J = \{\pm 1, \pm 3, \ldots, \pm j_{\max}\}$ for some odd $ j_{\max}\in\IN $
in \eqref{Ansatz.discussion} and \eqref{PDE.mfe}. We conjecture that it is indeed possible to prove that 
for $ j_{\max}=5 $ the error is $\Ord{\epsi^3}$ on $ [0,\tstar/\epsi] $ if our assumptions are adapted in an obvious way. The new terms would make the proofs even more complicated, but the techniques and the strategy would not have to be changed significantly. 
In addition to 
\eqref{Bounds.uhat} and \eqref{Bounds.Dmu.uhat} one has to show that 
\begin{align*}
\sup_{t\in[0,\tstar/\epsi]}\| \uhat_5(t) \|_{L^1} 
\leq C \epsi^2
\qquad \text{and} \qquad
\sup_{t\in[0,\tstar/\epsi]}\| D_\mu \uhat_5(t) \|_{L^1} 
\leq C \epsi^2
\end{align*}
in the framework of Propositions~\ref{Proposition01} and \ref{Proposition02}.
We believe, however, that unfortunately an extension to $ j_{\max}=7 $ and an error of $\Ord{\epsi^4}$ is \emph{not} possible with our techniques. The reason is, roughly speaking, that for $ j_{\max}=7 $ the difficult case in 
Step~4 of the proof of Theorem~\ref{Theorem.error.bound} is $ |J|_1 = |j| = 9 $. This case appears, in particular, for $J=(3,3,3)$, and in order to handle this case, we would need a non-resonance condition which is not fulfilled for the Maxwell--Lorentz system and the Klein--Gordon system.

\bibliographystyle{abbrv}
\bibliography{hfwp}

\end{document}